\documentclass[11pt, a4paper]{article}
\usepackage[english]{babel}

\usepackage{amsmath}
\usepackage{amssymb}
\usepackage{amsthm}
\usepackage{mathabx}
\usepackage{tikz}
\usepackage{tikz-cd}
\usepackage{todonotes}
\usepackage{hyperref}
\usepackage{cleveref}
\usepackage{footmisc}
\usepackage{bbm}
\usepackage{afterpage}
\usepackage{mathrsfs}
\usepackage{sseq}
\usepackage{spectralsequences}
\usepackage{extarrows}
\usepackage{enumitem}
\usepackage{stmaryrd}
\usepackage{pdflscape}

\newcommand{\gr}{\mathrm{gr}}

\newcommand{\ord}{\mathrm{ord}}

\newcommand{\1}{\mathbf{1}}
\newcommand{\BP}{\mathrm{BP}}

\renewcommand{\L}{\mathrm{L}}
\DeclareMathOperator{\KU}{KU}
\DeclareMathOperator{\ku}{ku}
\DeclareMathOperator{\KO}{KO}

\DeclareMathOperator{\ko}{ko}

\DeclareMathOperator{\MU}{MU}

\newcommand{\Sph}{\mathbf{S}}

\DeclareMathOperator{\tmf}{tmf}
\DeclareMathOperator{\Tmf}{Tmf}

\DeclareMathOperator{\Mod}{Mod}

\DeclareMathOperator{\Sp}{Sp}

\DeclareMathOperator{\Syn}{Syn}

\newcommand{\colim}{\mathrm{colim}\,}
\DeclareMathOperator{\Ext}{Ext}

\newcommand{\id}{\mathrm{id}}

\DeclareMathOperator{\map}{map}

\newcommand{\A}{\mathcal{A}}

\newcommand{\C}{\mathcal{C}}

\newcommand{\E}{\mathbf{E}}

\newcommand{\F}{\mathbf{F}}

\newcommand{\h}{\mathrm{h}}

\renewcommand{\j}{\mathrm{j}}
\newcommand{\J}{\mathrm{J}}
\newcommand{\K}{\mathrm{K}}

\newcommand{\R}{\mathbf{R}}

\newcommand{\T}{\mathrm{T}}

\newcommand{\Z}{\mathbf{Z}}

\newcommand{\al}{\alpha}
\newcommand{\be}{\beta}

\usepackage{mathtools}

\theoremstyle{theorem}\numberwithin{equation}{section}
\newtheorem{theorem}[equation]{Theorem}
\crefname{theorem}{{th}.\!\!}{{ths}.\!\!}
\Crefname{theorem}{{Th}.\!\!}{{Ths}.\!\!}
\newtheorem{theoremalph}{Theorem}

\crefname{theoremalph}{{th}.\!\!}{{ths}.\!\!}
\Crefname{theoremalph}{{Th}.\!\!}{{Ths}.\!\!}

\Crefname{problem}{{Prb}.\!\!}{{Prbs}.\!\!}
\newtheorem{prop}[equation]{Proposition}
\Crefname{prop}{{Pr}.\!\!}{{Prs}.\!\!}
\newtheorem{lemma}[equation]{Lemma}
\Crefname{lemma}{{Lm}.\!\!}{{Lms}.\!\!}
\newtheorem{cor}[equation]{Corollary}
\Crefname{cor}{{Cor}.\!\!}{{Cors}.\!\!}

\Crefname{conjecture}{{Conj}.\!\!}{{Conjs}.\!\!}

\theoremstyle{definition}\numberwithin{equation}{section}
\newtheorem{mydef}[equation]{Definition}
\Crefname{mydef}{{Df}.\!\!}{{Dfs}.\!\!}

\Crefname{recall}{{Rcl}.\!\!}{{Rcls}.\!\!}

\Crefname{construction}{{Con}.\!\!}{{Cons}.\!\!}

\Crefname{assumption}{{As}.\!\!}{{As}.\!\!}

\Crefname{notation}{{Nt}.\!\!}{{Nts}.\!\!}

\Crefname{situation}{{St}.\!\!}{{Sts}.\!\!}

\theoremstyle{remark}\numberwithin{equation}{section}
\newtheorem{example}[equation]{Example}
\Crefname{example}{{Ex}.\!\!}{{Exs}.\!\!}

\Crefname{nonexample}{{NonEx}.\!\!}{{NonEx}.\!\!}

\Crefname{claim}{{Clm}.\!\!}{{Clms}.\!\!}
\newtheorem{remark}[equation]{Remark}
\Crefname{remark}{{Rmk}.\!\!}{{Rmks}.\!\!}

\Crefname{idea}{{Id}.\!\!}{{Ids}.\!\!}
\newtheorem{warn}[equation]{Warning}
\Crefname{warn}{{Warn}.\!\!}{{Warns}.\!\!}

\Crefname{figure}{{Fig.}\!\!}{{Figs.}\!\!}
\Crefname{footnote}{{Fn.}\!\!}{{Fn.}\!\!}

\Crefname{part}{{\textsection}\!\!}{{\textsection}\!\!}
\Crefname{chapter}{{\textsection}\!\!}{{\textsection}\!\!}
\Crefname{section}{{\textsection}\!\!}{{\textsection}\!\!}
\Crefname{subsection}{{\textsection}\!\!}{{\textsection}\!\!}
\Crefname{appendix}{{\textsection}\!\!}{{\textsection}\!\!}

\usepackage{dutchcal}

\begin{document}
\title{A synthetic approach to detecting $v_1$-periodic families}
\author{
Christian Carrick\footnote{\url{c.d.carrick@uu.nl}}\,\, and Jack Morgan Davies\footnote{\url{davies@math.uni-bonn.de}}
}
\maketitle

\begin{abstract}
We provide a new proof that the unit map from the sphere spectrum to the connective image-of-$J$ spectrum $\j$ is surjective on homotopy groups. This is done by constructing modified Adams and Adams--Novikov spectral sequences for $\j$ using synthetic spectra. We easily compute these spectral sequences in their entirety using long exact sequences and Adams operations on topological $K$-theory. In particular, these computations and deductions are carried out without ever directly calculating $\F_p$- or $\BP$-homology nor the associated Ext groups of $\j$ nor $\ko$.
\end{abstract}

\setcounter{tocdepth}{2}
\tableofcontents

\section{Introduction}

For each prime $p$, let $\j$ denote the \emph{connective image-of-$J$ spectrum}, defined by the fibre sequence of spectra
\[\j \to \ko \xrightarrow{\psi^3-1} \tau_{\geq 4} \ko\]
at the prime $2$, for example. Using synthetic spectra, we construct modified Adams and Adams--Novikov spectral sequences for $\j$ and compute them entirely---which is a simple task, on par with computing the homotopy groups of $\j$. We do this without ever computing the $\F_p$- or $\BP$-homology of $\j$ nor the associated $\Ext$-groups. This is a drastic simplification of previous methods to calculate related spectral sequences for $\j$. Indeed, the classical Adams spectral sequence (ASS) for $\j$ requires difficult and subtle computations, going back to Davis \cite{daviscohomologyofj} and completed only recently by Bruner--Rognes \cite{imageofjbrunerrognes}. As far as the authors are aware, there is no computation of an Adams--Novikov spectral sequence (ANSS) of $\j$ in the literature nor of the comodule $\BP_\ast \j$.\\

Of course, the homotopy groups of $\j$ are well-known and easy to compute from the homotopy groups of $\ko$ and the action of the Adams operations. The utility of these modified spectral sequences is not as another approach to compute the homotopy groups of $\j$, but rather as a tool to prove detection statements. To demonstrate this, we provide a short and simplified proof of the following classical fact, widely credited to Adams \cite{adamsjx} and unpublished work of Sullivan (see May \cite[\textsection17]{infiniteloopspaces} and Madsen \cite[Th.5.5]{ibmadsenjstatement}).

\begin{theoremalph}\label{maintheorem}
Fix a prime $p$. The unit map $\Sph\to \j$ induces a surjection on homotopy groups.
\end{theoremalph}

The existence of the Adams $e$-invariant shows that the surjection above is in fact split. Our modified Adams--Novikov spectral sequence for $\j$ used to prove \Cref{maintheorem} can also be used to refine this splitting to one of spectral sequences (\Cref{splitonlatterpages}) at odd primes and a split surjection of filtered abelian groups (\Cref{filtereredsplitness}) at $p=2$. One of our key tools is a novel $t$-structure on the $\infty$-category of synthetic spectra (\Cref{abstractmaintheorem}), which we will continue to explore in future work.
\\

The background and basic facts about synthetic spectra and synthetic forms of topological $K$-theory take up a majority of this article. If one takes the Adams and Adams--Novikov spectral sequences for $\ko$ for granted, then our short proofs of \Cref{maintheorem,splitonlatterpages,filtereredsplitness} are mostly concentrated in \Cref{imageofjconstructionsection,selfmapssection,detectionresultsforj}.

\subsection{History and motivation}

The \emph{$J$-homomorphism} is a nontrivial homomorphism from the homotopy groups of the infinite special orthogonal group to the stable homotopy groups of spheres $\pi_\ast\Sph$. In \cite{adamsjx}, Adams uses real topological $K$-theory $\KO$ to study the $J$-homomorphism, and proved two facts about $\pi_\ast\Sph$ at the prime $2$:

\begin{itemize}
\item The image of the $J$-homomorphism detects a certain cyclic subgroup of $\pi_k\Sph$ for $k\geq 1$ congruent to $0,1,3$, or $7$ modulo $8$. The order of this subgroup was determined by Adams in \cite{adamsjx} up to a factor of $2$ in degrees congruent to $7$ modulo $8$; this ambiguity is known as the \emph{Adams conjecture} and was confirmed by Quillen \cite{quillenadams}.
\item Writing $\ko$ for the \emph{connective real $K$-theory} spectrum, the unit map $\Sph\to \ko$, arising from the fact that $\ko$ is a ring spectrum, detects certain cyclic subgroups of order $2$ inside $\pi_k\Sph$ for $k\geq 1$ congruent to $1$ or $2$ modulo $8$. These subgroups are generated by the so-called \emph{$\mu$-family}.
\end{itemize}

Let us see in detail how one can prove that $\ko$ detects the $\mu$-family $\al_{1+4k} \in \pi_{8k+1} \Sph$ at the prime $2$. First, recall that Adams constructs a self-map $v_1^4 \colon \Sigma^8 \Sph/2 \to \Sph/2$ of the mod $2$ Moore spectrum which induces multiplication by $\overline{u}^4$ on $\KU$-homology
\[\pi_\ast (v_1^4 \otimes \KU) \colon \F_2[u^\pm] \simeq \pi_\ast \KU/2 \to \pi_{\ast+8} \KU/2,\]
where $\overline{u}$ is the generator of $\pi_2 \KU/2$. This can be done by lifting $8\sigma \in \pi_7 \Sph$ through the boundary map
\[\partial \colon \pi_8 \Sph/2 \to \pi_7 \Sph,\]
and that $\pi_8 \Sph/2$ is $2$-torsion, so that this lift $\Sigma^8\Sph \to \Sph/2$ factors as our desired $v_1^4$. Writing $\widetilde{\eta}\colon \Sigma^2\Sph \to \Sph/2$ for a lift of $\eta$ through the boundary map $\partial\colon\pi_2 \Sph/2 \to \pi_1 \Sph$, we then define
\[\al_{1+4k} \colon \Sph^{8k+1} \xrightarrow{\widetilde{\eta}} \Sigma^{8k-1} \Sph/2 \xrightarrow{(v_1^{4})^k} \Sigma^{-1} \Sph/2 \xrightarrow{\partial} \Sph\]
as an element in $\pi_{8k+1} \Sph$. It is not a priori clear of course, that $\al_{1+4k}$ is nonzero. Adams proves that this class is nonzero using his $e$-invariant and Toda brackets; see \cite[Th.12.13]{adamsjx}.\\

Another attempt to prove that these classes are nonzero is as follows: first, we notice that $x_k=(v_1^4)^k \circ \widetilde{\eta}$ is sent to a nonzero element $\overline{u}^{4k+1}$ in $\pi_{8k+2} \KU/2$, essentially by construction. As the sphere $\Sph$ is connective, it suffices to work with the connective form of $K$-theory $\ku = \tau_{\geq 0}\KU$. The unit $\Sph \to \ku$ factors through $\ko$, so these classes $x_k$ must also be nonzero in $\pi_{8k+2} \ko/2$. The diagram
\[\begin{tikzcd}
    {\Sph}\ar[r, "{\partial_\Sph}"]\ar[d, "\overline{h}"]  &   {\Sigma\Sph/2}\ar[d, "h"]    \\
    {\ko}\ar[r, "{\partial_{\ko}}"]         &   {\Sigma \ko/2}
\end{tikzcd}\]
commutes, so to see that $\al_{1+4k}$ is nonzero, it suffices to show that $\partial_{\ko}(\overline{h}(x_k))$ is nonzero in $\pi_{8k+1}\ko$. This is not immediately clear though, as in the exact sequence
\[\pi_{8k+2} \ko \xrightarrow{q} \pi_{8k+2} \ko/2 \xrightarrow{\partial_{\ko}} \pi_{8k+1} \ko\]
one must prove that the class $\overline{h}(x_k)$ does not lie in the image of $q$.\\

This line of reasoning can work more smoothly if one argues using long exact sequences between ASSs or ANSSs, as the extra filtration can help to separate the image of $q$ and $\overline{h}(x_k)$. The appearance of spectral sequences is also quite natural: one way to prove that the map $v_1^4$ induces an isomorphism on $\KU$-homology is to use the Adams--Novikov spectral sequence. This approach has two key difficulties though: 
\begin{itemize}
    \item one needs to compute the $E_2$-pages of the associated spectral sequences (perhaps even the whole spectral sequences), and
    \item one needs to come to terms with the subtle and sometimes limited \emph{geometric boundary theorem} that relates the boundary maps above to algebraic boundary maps on $E_2$-pages; see \cite[Th.2.3.4]{greenbook}.
\end{itemize}
These problems are manageable in the above example for $\ko$. However, $\ko$ cannot be used to detect all the classes in the image-of-$J$, at the very least as $\pi_d \ko$ vanishes for $d\equiv 3,7$ modulo $8$.\\

A spectrum that does end up detecting both the image-of-$J$ and the $\mu$-family is the aptly named \emph{connective image-of-$J$} spectrum, simply denoted by $\j$. This plays a key role in Mahowald's study of $v_1$-periodic phenomena in \cite{mahowaldJinEHP}. The fact that $\j$ detects these two families is essentially the content of \Cref{maintheorem}. One definition of $\j$ at the prime $2$ is via the the fibre sequence
\begin{equation}\label{fibreseq}\j\to \ko\xrightarrow{\psi^3-1} \tau_{\geq 4}\ko,\end{equation}
where $\psi^3$ is the third stable Adams operation---the $4$-connective cover $\tau_{\geq 4}\ko$ can be seen as making up for the fact that the image of the $J$-homomorphism and the $\mu$-family agree in $\pi_1$ and $\pi_2$. The homotopy groups of $\ko$ are well-known by Bott periodicity, as is the action of the Adams operations $\psi^k$ on $\pi_\ast\ko$, so the homotopy groups of $\j$ are readily calculated from the above fibre sequence.\\

Unfortunately, the problems with the spectral sequence approach to detection outlined above apply even more drastically to $\j$. The ASS for $\j$ is much more complicated than that of $\ko$, with its $E_2$-page computed by Davis \cite{daviscohomologyofj} and the rest of the structure only obtained recently by Bruner--Rognes \cite{imageofjbrunerrognes}. As far as the authors are aware, there is no record of a computation of the comodule $\BP_\ast \j$ in the literature, let alone a computation of the ANSS of $\j$.\\

This is where \emph{synthetic spectra} come in, as a method for producing \emph{modified} ASSs and ANSSs. In particular, modified spectral sequences for $\j$. Once appropriate modified spectral sequences are in place for $\j$, designed to be easily computable, the detection arguments outlined above run through seamlessly to prove \Cref{maintheorem}.

\subsection{Synthetic techniques}
The $\infty$-category of synthetic spectra $\Syn_E$ is a categorification of the $E$-based Adams spectral sequence ($E$-ASS) for some Adams-type homology theory $E$. Each synthetic spectrum $X$ has an associated spectral sequence $\sigma X$. In particular, for each spectrum $X$ there is a synthetic spectrum $\nu X$ whose associated spectral sequence is the standard $E$-ASS for $X$. There can be, however, many other synthetic spectra $Y$ not of the form $\nu X$ whose associated spectral sequence still abuts to the homotopy groups of $X$---we call such a spectral sequence a \emph{modified $E$-ASS} for $X$; see \Cref{signaturedefinition}.\\

To produce a modified $E$-ASS for $\j$, we would like to repeat the definition of $\j$ using the fibre sequence (\ref{fibreseq}). To do that, we introduce a new $t$-structure on $\Syn_E$; see \Cref{tstructuresexist} for a more general statement.

\begin{theoremalph}\label{abstractmaintheorem}
Let $E$ be a connective $\E_1$-ring of Adams-type. There exists a left and right complete monoidal $t$-structure on $\Syn_E$, called the \emph{vertical $t$-structure}, whose connective objects are given by those synthetic spectra $X$ with $\pi_{a,b}X=0$ for $a<0$. Moreover, the heart of this $t$-structure is given by the abelian category of graded $\pi_{0,\ast}\1$-modules.
\end{theoremalph}

We then define $\j_E$ using the fibre sequence
\[\j_E \to \nu \ko \xrightarrow{\psi^3-1} \tau_{\geq 4}^\uparrow \nu \ko\]
for both $E=\F_2$ and $\BP$ at the prime $2$, where $\tau_{\geq 4}^\uparrow$ refers to the $4$-connective cover with respect to this vertical $t$-structure. As the ASS and ANSS for $\ko$ are well-known, the spectral sequence to $\j_E$ can be easily deduced from the long exact sequence on mod $\tau$ bigraded homotopy groups induced by the above fibre sequence. In fact, we argue that the ASS and ANSS for $\ku$ suffice: the ASS and ANSS for $\ko$ also follow from purely synthetic methods and Wood's theorem. In particular, the knowing Wood's theorem and the very simple ASS and ANSS of $\ku$, one can compute the spectral sequence associated to $\j_E$ in its entirety. The same also holds at odd primes.\\

Aside from just proving \Cref{maintheorem}, synthetic methods also give refinements. The Adams $e$-invariant shows that $\pi_\ast \Sph \to \pi_\ast \j$ is split surjective in each degree. For odd primes $p$, this generalises to a split surjection of (modified) Adams--Novikov spectral sequences.

\begin{theoremalph}\label{splitonlatterpages}
For odd primes $p$, the unit map $\1\to \j_\BP$ induces a split surjection on synthetic homotopy groups, and at $p=2$, it induces a split surjection on $E_\infty$-pages of the associated spectral sequences.
\end{theoremalph}

Although it is not true that $\1\to \j_\BP$ induces a surjection on bigraded homotopy groups at the prime $2$, see \Cref{failiureofsurjectivity}, the modified ANSS associated with $\j_{\BP}$ does allow us to strengthen the classical statement that $\pi_\ast\Sph\to \pi_\ast\j$ is a split surjection of abelian groups.

\begin{theoremalph}\label{filtereredsplitness}
Fix $p=2$. Consider $\pi_\ast\Sph$ and $\pi_\ast\j$ as filtered abelian groups using the (modified) ANSSs associated with the $\BP$-synthetic spectra $\1$ and $\j_\BP$. Then the map $\pi_\ast \Sph\to \pi_\ast \j$ is a split surjection of \emph{filtered} abelian groups.
\end{theoremalph}

As a corollary of our computation of $\pi_{\ast,\ast}\j_{\F_2}$, we immediately obtain a synthetic analogue of a result of Andrews--Miller \cite{invertingeta} and Miller--Ravenel--Wilson \cite{mrw} that the synthetic $\al$-family $\al_n$ supports $\al_1^k$-multiplication for all $k\geq 0$ and $n\neq 2$; see \Cref{millerandrewsthm}.\\

We also use synthetic spectra to explore the periodic spectra $\KO$ and $\J$, again with an eye to producing computable and useful modified ASSs and ANSSs. This is particularly interesting for the ASS, as the classical ASSs for all of the periodic spectra $\KU$, $\KO$, and $\J$ are trivial and hence do not converge. This modified ASS for $\J$ behaves like a nontrivial ASS converging to the $\K(1)$-local sphere; see \Cref{ktheoryconstructionsection} and \Cref{periodicvariantsofJsubsection}.

\subsection*{Outline}

This article is divided into three main sections: in \Cref{taugabsssection,tstructuresection}, we discuss the $\infty$-category $\Syn$ and provide various new $t$-structures, then in \Cref{ktheoryconstructionsection,imageofjconstructionsection}, we explore synthetic lifts of topological $K$-theory and image-of-$J$ spectra providing us with the needed modified $E$-ASSs, and finally in \Cref{selfmapssection,detectionresultsforj}, we prove the desired detection results using these synthetic constructions. Experts on synthetic spectra may comfortably skip \Cref{taugabsssection}, and experts on topological $K$-theory and the Hurwicz image of $\ko$ may skip \Cref{ktheoryconstructionsection,selfmapssection}---these latter two sections serve to provide motivation and notation for our study of $\j$. The reader who wants a simple proof of \Cref{maintheorem} can focus on sections \Cref{imageofjconstructionsection,selfmapssection,detectionresultsforj}. In some more detail:

\begin{itemize}
\item In \Cref{taugabsssection}, the connection between synthetic spectra, filtered spectra, and their associated spectral sequences is discussed. In particular, this provides us with a language to discuss \emph{modified} $E$-ASSs, reproducing $v_n$-localised ASSs as an example.
\item In \Cref{tstructuresection}, we construct a collection of linear $t$-structures on synthetic spectra and prove generalisations of \Cref{abstractmaintheorem}---\emph{linear}, as they are based on cutting synthetic spectra using lines in $\mathbf{R}^2$. The heart of these $t$-structures is calculated and we prove other basic properities of the associated truncation and connective cover functors.
\item In \Cref{ktheoryconstructionsection}, various $E$-synthetic lifts of topological $K$-theory are constructed, first for $E=\F_p$ and then for $E=\BP$. First, we use synthetic spectra to prove Wood's theorem from first principles. Then we show how the classical ASSs and ANSSs for $\KO$ and $\ko$ can be recovered from a synthetic version of Wood's theorem using only simple long exact sequence arguments. These constructions are a warm-up to our synthetic construction of various image-of-$J$ spectra.
\item In \Cref{imageofjconstructionsection}, the main characters of this article are introduced: the synthetic lifts for $\j$. Similar to the classical calculation of the homotopy groups $\j$, the synthetic homotopy groups of these lifts are then easily computed simultaneously giving us modified ASSs and ANSSs for $\j$. We also show that each $\j_E$ has a chosen $\E_\infty$-structure in $\Syn$ and produce periodic variants.
\item In \Cref{selfmapssection}, we reprove the classical statements of Adams concerning the existence of $v_1$ self-maps on mod $p$ Moore spectra and how these can be used to prove detection results concerning real topological $K$-theory. These arguments help set the stage for following section.
\item In \Cref{detectionresultsforj}, we prove \Cref{maintheorem}, stating the surjectivity of the unit map $\Sph\to \j$ on homotopy groups, using the calculations from \Cref{imageofjconstructionsection} together with the detection arguments for $\ko$ from \Cref{selfmapssection}. We also prove two more refined statements: \Cref{splitonlatterpages}, which states that the synthetic unit map $\1\to \j_{\BP}$ induces a \textbf{split} surjection on synthetic homotopy groups at odd primes, and \Cref{filtereredsplitness}, which states that the surjection of \Cref{maintheorem} refines is a split surjection of filtered abelian groups at the prime $2$.
\item In \Cref{sssection}, we provide charts for the various modified $E$-ASS used throughout this article.

\end{itemize}

\subsection*{Other work}

The study of $p$-complete even $\MU$-synthetic spectra through the guise of \emph{cellular $\mathbf{C}$-motivic stable homotopy theory} has already led others to study objects closely related those appearing in this paper; the connection between motivic homotopy theory and $\MU$-synthetic spectra appears in Pstragowski's original article \cite[Th.1.4]{syntheticspectra}. For example, in \cite{isasksenshkembi}, Isaksen--Shkembi study both connective covers for motivic spectra in \textsection 3 as well as motivic homotopy groups of motivic variants of $K$-theory spectra, both periodic and connective in \textsection 4---unsurprisingly, our results concerning the synthetic homotopy groups of various synthetic $K$-theory spectra look the same. \\

There has also recently been an explosion in work related to synthetic spectra and filtered spectra. In particular, the well-known, but not yet documented, result about the cellularity of $\F_p$-synthetic spectra proven in \Cref{cellularityoffp} has recently been generalised by Lawson \cite{lawson2024synthetic}. Indeed, he shows that the $\infty$-category of $E$-synthetic spectra is cellular for all connective ring spectra $E$ of Adams type. Moreover, Lee--Levi \cite{leelevi} also construct certain $t$-structures on filtered objects in a general stable $\infty$-category $\C$ similar to our linear $t$-structure of \Cref{tstructuresection}.\\

The general philosophy of working with modified $E$-ASSs in synthetic spectra has been a theme in many paper, such as \cite{burkhahnseng,burklundmoore}, just to mention a few. In \cite{smfcomputation}, for example, together with van Nigtevecht we compute the associated spectral sequence, there called the \emph{signature}, of the synthetic spectrum $\mathrm{Smf}$, a synthetic lift of the spectrum of projective topological modular forms $\Tmf$, and state that it is not clear how to calculate the classical ANSS for $\Tmf$ from first principles; see \cite[\textsection 1.2.1]{smfcomputation}. The identification of various modified $E$-ASSs, or signatures, is also further explored in \cite[\textsection1-2]{osyn}. 

\subsection*{Future work}

To reiterate, our proof of \Cref{maintheorem} does not use any homology calculations of $\j$. One could attempt to reprove the results of this paper using the classical ASS for $\j$ as computed in \cite{imageofjbrunerrognes}, although this would have rely on $v_1$-periodic techniques in the ASS, as opposed to the direct approach we take with a modified ANSS for $\j$. There are, however, other situations where direct ASS computations become unmanageable. The authors' initial motivation for these synthetic constructions was to try to calculate the image of the unit map from $\Sph$ into the fibre of a map $F-\lambda\colon\tmf\to \tau_{\geq 8}\tmf$, where $F$ is an Adams operation $\psi^k$ or a Hecke operator $\T_n$, from either \cite{adamsontmf}, \cite{realspectra}, or \cite{heckeontmf}, and $\lambda$ is the eigenvalue of the unit. Computing directly with the ASS for $\tmf$ is very difficult, as is determining the homology of $\tau_{\geq 8}\tmf$, let alone its ASS. We have had success using similar techniques as explored in this article to compute the Hurewicz image of the fibre of $\psi^2-1 \colon \tmf \to \tmf$ at the prime $3$; see \cite{heighttwojat3}. We are currently exploring the fibre of $\psi^3-1$ on $\tmf_2$ at the prime $2$ as well as computations related to Behrens' $\mathrm{Q}(N)$-spectra, which are another variant of image-of-$J$ spectra at height $2$; see \cite{ktwospheremark}.

\subsection{Notation}
Throughout this article, $p$ is a prime and our spectra are implicitly $p$-complete. If $E$ is an Adams-type homology theory, we write $E$-ASS for the $E$-based Adams spectral sequence. We further abreviate $\F_p$-ASS to ASS, and $\BP$-ASS to ANSS. We denote the category of $E$-synthetic spectra $\Syn_E$ of \cite{syntheticspectra} be $\Syn$. We use the ``stem--filtration" grading for synthetic spectra, meaning that $\Sigma^{s,f}\1=\Sigma^{-f}\nu\Sph^{s+f}$. We will also implicitly use the identification of groups
\[\pi_{s,f} \nu X /\tau \simeq \Ext_{E_\ast E}^{f,s-f}(E_\ast X, E_\ast)\]
for a spectrum $X$, due to Pstragowski \cite[Lm.4.56]{syntheticspectra}; also see \cite[Pr.1.25]{osyn} and \Cref{examplesofsignatures}.

\subsection{Basic assumptions}\label{sec:basicassumptions}
We base this article on minimal assumptions about complex topological $K$-theory and the stable homotopy groups of spheres. More specifically, we will use the following simple facts:

\begin{itemize}
    \item There is a $\E_\infty$-ring $\KU$ of \emph{complex topological $K$-theory} such that $\pi_\ast \KU \simeq \Z[u^\pm]$ with $|u|=2$, and as a complex-oriented cohomology theory, $\KU$ is equipped with the multiplicative form group law.
    \item There is an action of $\Z_p^\times$ on this $\E_\infty$-ring $\KU$ such that for each $\lambda \in \Z_p^\times$, the action of $\lambda$ on $\pi_2\KU$ is multiplication by $\lambda$. We write the associated endomorphism of $\KU$ as $\psi^\lambda$.
    \item The computations of the first $8$ stable homotopy groups of spheres, so $\pi_d \Sph$ for $d\leq 8$, as well as the ASS and ANSS for $\Sph$ in this range. Moreover, we assume the computation of $\pi_2 \Sph/2 \simeq \Z/4\Z$.
\end{itemize}

As we are implicitly $p$-completing everywhere, one can construct such a $\KU$ with such an $\Z_p^\times$-action as the Morava $E$-theory $E_1$ associated to the formal multiplicative group over $\F_p$; see \cite{gh04}, \cite[\textsection7]{pp}, or \cite[Con.5.1.1]{ec2name} for such constructions. The stable homotopy groups of spheres in this range follow from straight forward computation; see \cite{novicesanss}. The computation of $\pi_2 \Sph/2$ follows from an explicit ASS computation, for example.\\

As a consequence of these assumptions, we have the following definition.

\begin{mydef}
    Let $\KO=\KU^{hC_2}$, where $C_2$ is the subgroup of $\Z_p^\times$ generated by $-1$. At odd primes, let $\L = \KU^{h\F_p^\times}$ where $\F_p^\times$ is the maximal finite subgroup of $\Z_p^\times$. We write $\ku = \tau_{\geq 0}\KU$, $\ko=\tau_{\geq0} \KO$, and $\ell = \tau_{\geq 0} \L$, which all inherit Adams operations by functoriality.
\end{mydef}

\begin{remark}\label{rmk:easyconsequences}
The following elementary facts about the above spectra will also be used.

\begin{enumerate}
    \item The ANSS for $\ku$ has $E_2$-page $\Z[u]$ with $|u|=(2,0)$ and collapses immediately with no differentials. Similarly, at odd primes, the ANSS for $\ell$ has $E_2$-page $\Z[v_1]$ with $|v_1|=2p-2$ and collapses immediately with no differentials. This follows from the facts for $E=\ku$ or $\ell$, that are complex-oriented. In this case, $E$ is a homotopy $\BP$-module and the associated coaugmented cosimplicial AN-tower $E\otimes \BP^{\otimes(\bullet+1)}$ has an extra degeneracy.
    \item There are isomorphisms of $\A_\ast$-comodules
    \[H_\ast(\ku;\F_2) \simeq \A_\ast \square_{\mathcal{E}(1)_\ast} \F_2, \quad H_\ast(\ell;\F_p) \simeq \A_\ast \square_{\mathcal{E}(1)_\ast} \F_p\]
    at $p=2$ and odd primes $p$, respectively, where $\mathcal{E}(1)_\ast$ is the dual to the subalgebra of $\A$ generated by the Milnor elements $Q^0$ and $Q^1$. This follows from standard calculations of Wilson \cite[Pr.1.7]{wilsonbpii} using the cofibre sequences of spectra
\[\Sigma^{2p-2}\BP\langle 1\rangle\xrightarrow{v_1} \BP\langle 1\rangle\to \Z\qquad \qquad \Z\xrightarrow{p} \Z\to \F_p.\]
\end{enumerate}
\end{remark}





\subsection*{Acknowledgements}

We would like to thank Shaul Barkan, Emma Brink, Lennart Meier, Sven van Nigtevecht, Lucas Piessevaux, and John Rognes for their helpful conversations and suggestions. In particular, thank you to Shaul and Sven for catching some inaccuracies in some of our previous computations. Thank you to the anonymous referee for their helpful suggestions.\\

The first author was supported by NSF grant \texttt{DMS-2401918} as well as the NWO grant \texttt{VI.Vidi.193.111}. The second author is an associate member of the Hausdorff Center for Mathematics at the University of Bonn (\texttt{DFG GZ 2047/1}, project ID \texttt{390685813}). Both authors would like to thank the Isaac Newton Institute for Mathematical Sciences, Cambridge, for support and hospitality during the programme \emph{Equivariant homotopy theory in context} where the revision of this paper was undertaken. This work was supported by EPSRC grant no EP/K032208/1.


\section{The spectral sequence associated to a synthetic spectrum}\label{taugabsssection}

Associated with a synthetic spectrum $X$ are two spectral sequences: the \emph{$\tau$-Bockstein spectral sequence} (BSS) and what we call the \emph{$\sigma$-spectral sequence} ($\sigma$-SS). These concepts will be used to concretely define what we mean by \emph{modified $E$-ASSs}. Recall that we are using the ``stem--filtration" grading for synthetic spectra, meaning that $\Sigma^{s,f}\Sph=\Sigma^{-f}\nu\Sph^{s+f}$. In particular, $\tau$ has degree $(0,-1)$ and suspension as a stable $\infty$-category has degree $(1,-1)$.\\

This section uses filtered spectra as a homotopical model for spectral sequences as done in \cite[\textsection1.2.2]{haname}, \cite[\textsection II.1]{alicethesis}, and \cite[\textsection2]{mmf}, for example. Other references with an eye towards synthetic spectra can also be found in \cite[\textsection1.2]{osyn} and \cite[\textsection2.1]{smfcomputation}.


\subsection{Filtered spectra and spectral sequences} 

Spectral sequences come from filtrations. In stable homotopy theory, this is encoded as saying that spectral sequences come from filtered spectra. A filtered spectrum is an object in $\mathrm{Fil}(\Sp):=\mathrm{Fun}(\Z^{\mathrm{op}},\Sp)$, the $\infty$-category of filtered spectra, where $\Z$ is the poset consisting of the integers with the usual ordering. Thus a filtered spectrum consists of a tower $X^\bullet$ of spectra
\[\cdots \to X^{n+1}\to X^{n}\to X^{n-1}\to\cdots.\]
where $X^i\in \Sp$. One extracts a spectral sequence from $X^\bullet$ by forming the exact couple
\[\begin{tikzcd}
\bigoplus\limits_{s,t}\pi_{t-s}X^{t+1}\arrow[rr]&&\bigoplus\limits_{s,t}\pi_{t-s}X^t\arrow[dl]\\
&\bigoplus\limits_{s,t}\pi_{t-s}X_t^t\arrow[ul, dashed]
\end{tikzcd}\]
giving a spectral sequence with signature
\begin{equation}\label{ssofafilteredguy}E_2^{t-s,s}=\pi_{t-s}X_t^t\implies \pi_{t-s}X^{-\infty}\end{equation}
and differentials $d_r$ of bidegree $(-1,r)$. Here we write $X_t^t:=\mathrm{cofib}(X^{t+1}\to X^t)$ and $X^{-\infty}=\colim X^n$, and the (descending) filtration on $\pi_{t-s}X^{-\infty}$ is defined by
\[F^s\pi_{t-s}X^{-\infty}=\mathrm{image}(\pi_{t-s}X^t\to\pi_{t-s}X^{-\infty})\]
for all $t\in\Z$. If $\lim X^n=X^\infty=0$, then the spectral sequences converge in the sense that the $F^s$ filtration is exhaustive and Hausdorff, and that $E_\infty\simeq F^s/F^{s+1}$; see \cite[\textsection5]{boardman}.\\

If one knows the homotopy groups $\pi_*X^n$ and the effect of the maps $\pi_*X^n\to\pi_*X^{n-1}$ for every $n$, then one knows the $E_2$ page and the behaviour of all the differentials via the exact couple. It is helpful then to package the homotopy groups $\pi_*X^n$ as $n$ varies into one bigraded group. 

\begin{mydef}
For $X^\bullet\in\mathrm{Fil}(\Sp)$, we define the bigraded homotopy group $\pi_{n,s}X^\bullet=\pi_nX^{n+s}$.
\end{mydef}

\begin{remark}
This indexing convention does not agree with many sources on bigraded homotopy groups. We have chosen them so that elements in $\pi_{n,s}X^\bullet$ correspond to classes in $\pi_nX^{-\infty}$ detected in filtration $s$ of the spectral sequence associated with $X^\bullet$.
\end{remark}

There are bigraded spheres in $\mathrm{Fil}(\Sp)$ that corepresent these bigraded homotopy groups; one simply sets 
\[\Sph^{n,s}(k):=\begin{cases}\Sph^n&k\le n+s\\0&k>n+s\end{cases}\]
where all the transition maps are either zero or the identity. The transition maps in a filtered spectrum $\pi_*X^n\to \pi_*X^{n-1}$ give $\pi_{*,*}X^\bullet$ the structure of a bigraded module over the bigraded ring $\Z[\tau]$ where $|\tau|=(0,-1)$. This is encoded in the category  $\mathrm{Fil}(\Sp)$ by the bigraded self-map $\tau:\Sph^{0,-1}\to \Sph^{0,0}$ of the unit $\1=\Sph^{0,0}$, defined in the obvious way. Since multiplication by $\tau$ in $\pi_{*,*}X^\bullet$ corresponds to the transition maps in the exact couple above, we see that the $\Z[\tau]$-module $\pi_{*,*}X^\bullet$ encodes the data of the spectral sequence associated with $X^\bullet$. This sentiment is captured in the following statement, which is proven by the standard arguments for exact couples; see \cite[Th.A.1]{burkhahnseng}, for example.

\begin{theorem}\label{filspthm}
We have the following facts about the relationship between the spectral sequence associated with $X^\bullet\in\mathrm{Fil}(\Sp)$ and the $\Z[\tau]$-module $\pi_{*,*}X^\bullet$.
\begin{enumerate}
\item There is a short exact sequence
\[0\to(\pi_{n,s}X^\bullet)/\tau\to E_2^{n,s}\xrightarrow{\partial}(\pi_{n-1,s+2}X^\bullet)[\tau]\to0\]
where $(\pi_{n-1,s+2}X^\bullet)[\tau]$ is the kernel of multiplication by $\tau$ on $\pi_{n-1,s+2}X^\bullet$.
\item The subgroup $(\pi_{n,s}X^\bullet)/\tau\simeq Z_\infty^{n,s}$ of $E_2^{n,s}$ consists of permanent cycles in stem $n$ and filtration $s$.
\item The subgroup $B_r^{n,s}\subset Z_\infty^{n,s}$ of classes that are hit by a $d_r$-differential is precisely the set of $x\in Z_\infty^{n,s}$ admitting a lift $\tilde{x}\in\pi_{n,s}X^\bullet$ such that $\tau^{r-1}x=0$.
\item Suppose the map $\partial:E_2^{n,s}\to(\pi_{n-1,s+2}X^\bullet)[\tau]$ sends $x\mapsto \tau^{r-2}y\neq0$ where $y$ is not divisible by $\tau$. Then $d_2,\ldots,d_{r-1}$ vanish on $x$ and $d_r(x)$ is represented by the image of $y$ along the map
\[\pi_{n-1,s+r}X^\bullet\twoheadrightarrow (\pi_{n-1,s+r}X^\bullet)/\tau\hookrightarrow E_2^{n-1,s+r}\]
Conversely, if $d_2,\ldots,d_{r-1}$ vanish on $x$ and $d_r(x)=y\neq0$, then $y$ admits a lift to $E_2$ that admits a lift $\tilde{y}$ to $\pi_{n-1,s+r}X$ such that $\partial(x)=\tau^{r-2}\tilde{y}$.
\end{enumerate}
\end{theorem}


\subsection{Synthetic spectra and modified Adams spectral sequences}
Let $E$ be an Adams-type spectrum, then for any $X\in\Sp$, one has the $E$-ASS of $X$
\[\Ext_{E_*E}(E_*,E_*X)\implies \pi_*(X^{\wedge}_E)\]
converging to the $E$-nilpotent completion of $X$. For some applications, this particular spectral sequence is ineffective and could be replaced by a variety of other spectral sequences with similar formal behaviour. These are called \textit{modified} $E$-ASSs for $X$. Such spectral sequences have a rich history; for example, one often considers localised ASSs, such as the \emph{$v_n$-localised ASS}
\[v_n^{-1}\Ext_{\mathcal{A}_*}(\F_p,H_*(Y))\implies \pi_*(v_n^{-1}Y)\]
for a spectrum $Y$ with a $v_n$ self-map; for example, see \cite[\textsection2.2]{akquigley}. This spectral sequence converges under certain conditions on $Y$, but it is not the same as the $\F_p$-ASS of $v_n^{-1}Y$, as $H_*(v_n^{-1}Y)=0$. Another common example of a modified ASS is the one used to compute the $X$ homology of a Moore spectrum $\Sph/(p^{i_0},v_1^{i_1},\ldots,v_n^{i_n})$ for some spectrum $X$, where one has a spectral sequence
\[\Ext_{\mathcal{A}_*}(\F_p,H(i_0,\ldots,i_n)\otimes H_*(X;\F_p))\implies X_*(\Sph/(p^{i_0},v_1^{i_1},\ldots,v_n^{i_n}))\]
where $H(i_0,\ldots,i_n)$ is a Moore object $\1/(q_0^{i_0},\ldots,q_n^{i_n})$ formed inside the derived category $\mathcal{D}(\mathrm{Comod}(\mathcal{A}_*))$; for example, see \cite[\textsection8]{mmmm}. This $E_2$-page is much easier to work with, but is different to the $E_2$-page of the $\F_p$-ASS of $X\otimes \Sph/(p^{i_0},v_1^{i_1},\ldots,v_n^{i_n})$ if any $i_j>1$.\\

Classically, the construction of these modified ASSs is quite complicated and \textit{ad hoc}. The category $\Syn$ of $E$-based synthetic spectra makes it easy to construct such spectral sequences and we also find it is convenient to think about objects in $\Syn$ in this way. Very roughly speaking, objects in $\Syn$ a pairs $(X,\mathcal{E})$, where $X$ is a spectrum, and $\mathcal{E}$ is a choice of modified $E$-based Adams spectral sequence for $X$. This point of view can be made into a rigorous description of $\Syn$ (up to technical issues in $\Syn$ like cellularity and $\tau$-completeness of the unit); see the \emph{filtered model} of \cite[Pr.C.22]{burkhahnsengartintate}. More categorically, $\Syn$ can be considered as a bundle of spectral sequences living over $\Sp$
\[
\begin{tikzcd}
\Syn\arrow[d,"\tau^{-1}"']\\
\Sp\arrow[u,bend right,"\nu"',dashed]
\end{tikzcd}
\]
where the synthetic analog functor $\nu$ gives a canonical section of this fibration, sending $X\mapsto (X,E$-ASS$(X))$ where $E$-$\mathrm{ASS}(X)$ is the standard $E$-ASS of $X$. From this perspective, the fibration $\tau^{-1}\colon \Syn\to\Sp$ is projection onto the first factor, and we devote the rest of the section to describing projection onto the second factor, namely how to extract a modified $E$-ASS for $\tau^{-1}X$ from a synthetic spectrum $X$. We refer the reader again to \cite[Appendix C]{burkhahnsengartintate} for a more categorical point of view on this construction.\\

The unit in $\Syn$ has a bigraded self-map $\tau:\Sigma^{0,-1}\1\to\1$, so for any $X\in\Syn$, one has a filtered synthetic spectrum
\[X^\bullet=\qquad \cdots\to\Sigma^{0,-1}X\xrightarrow{\tau} X\xrightarrow{\tau}\Sigma^{0,1}X\to\cdots\]
where $X^0=X$. Since $\Syn$ is a stable $\infty$-category, we may apply the mapping spectrum functor $\mathrm{map}_{\Syn}(\1,-)$ to this filtered synthetic spectrum to get a filtered spectrum we will call $\sigma(X)$.

\begin{theorem}\label{Gammathm}
The functor $\sigma:\Syn\to\mathrm{Fil}(\Sp)$ has the following properties
\begin{enumerate}
\item For any $X\in \Syn$, the map $\pi_{*,*}X\to\pi_{*,*}\sigma(X)$ is an isomorphism of $\Z[\tau]$-modules.
\item The spectral sequence associated with the filtered spectrum $\sigma(X)$ has signature
\[E_2=\pi_{t-s,s}(X/\tau)\implies \pi_{t-s}\tau^{-1}X.\]
This is the \emph{$\sigma$-SS} for $X$.
\item If $X$ is $\tau$-complete, the spectral sequence converges in the sense that the corresponding filtration on $\pi_*\tau^{-1}X$ is exhaustive and Hausdorff.
\item Considering $\mathrm{Fil}(\Sp)$ as a symmetric monoidal $\infty$-category with the Day convolution monoidal structure, then $\sigma$ has the canonical structure of a lax symmetric monoidal functor.
\end{enumerate}
\end{theorem}

\begin{proof}
The second and third claims follow from the first by definition of $\sigma$ along with the discussion in the previous section. For the first claim, one has
\begin{align*}
\pi_{n,s}\sigma(X)&=\pi_n\sigma(X)^{n+s}\\
&=\pi_n\mathrm{map}_{\Syn}(\1,\Sigma^{0,-n-s}X)\\
&=[\Sigma^{n,-n} \1,\Sigma^{0,-n-s}X]_{\Syn}\\
&=\pi_{n,s}X
\end{align*}
using that the stable $\infty$-categorical suspension functor on $\Syn$ coincides with the bigraded suspension $\Sigma^{1,-1}$. The lax monoidality of $\sigma$ comes from the natural maps
\[(\sigma X\otimes \sigma Y)^n\simeq \bigoplus_{i+j=n}\map_{\Syn}(\1,\Sigma^{0,-i}X)\otimes \map_{\Syn}(\1,\Sigma^{0,-j}Y)\]
\[ \to \bigoplus_{i+j=n}\map_{\Syn}(\1,\Sigma^{0,-i-j}X\otimes Y)\to \sigma(X\otimes Y)^n\]
and the unit map for $n\geq 0$
\[\Sph^{0,0}(n)=\Sph\to \map_{\Syn}(\Sigma^{0,n}\1,\1)\simeq \sigma(\1)^n.\qedhere\]
\end{proof} 

\begin{mydef}\label{signaturedefinition}
Given a synthetic spectrum $X$, we call the spectral sequence of \Cref{ssofafilteredguy} from the filtered spectrum $\sigma X$ the \emph{spectral sequence associated} to $X$, or the simply as the \emph{associated spectral sequence}, or the \emph{modified $E$-ASS}.
\end{mydef}

\begin{example}\label{examplesofsignatures}
When $Y=\nu X$, then \cite[Th.A.1]{burkhahnseng} shows that the spectral sequence associated to $Y$ is the standard $E$-ASS of $X$---this is also implicit in \cite{syntheticspectra}. Combining \Cref{filspthm,Gammathm}, we now have a precise meaning to the sense that the bigraded synthetic homotopy groups of $\nu X$ encode the data of the $E$-ASS of $X$. Another perspective on this fact can be found in \cite[Pr.1.25]{osyn}, where 
there is a natural equivalence of filtered spectra $\sigma(\nu X)$ and (the d\'{e}calage of) the $E$-ASS for $X$.
\end{example}

\begin{example}
    In \cite{osyn}, together with van Nigtevecht we construct various synthetic spectra $\mathcal{O}^{\mathrm{syn}}_{\mathsf{X}}(\mathsf{X})$ associated to certain spectral Deligne--Mumford stacks $(\mathsf{X}, \mathcal{O}^{\mathrm{top}}_{\mathsf{X}})$ whose associated spectral sequences are the \emph{descent spectral sequence}; see \cite[Th.B]{osyn}.
\end{example}

\begin{example}
Suppose $Y$ is a spectrum with a $v_n$ self-map $v$. Then using \cite[Th.A.1]{burkhahnseng}, one can show that there is a bigraded self-map $\tilde{v}$ of $\nu(Y)\in\Syn_{\F_p}$ such that $\tau^{-1}(\tilde{v})=v$ and reduction mod $\tau$ sends $\tilde{v}$ to the class in $\Ext_{\mathcal{A}_*}(\F_p,H_*Y)$ detecting $v$. It follows that the telescope $\tilde{v}^{-1}Y$ formed in $\Syn_{\F_p}$ has the property that
\[\pi_{*,*}(\tilde{v}^{-1}Y/\tau)=v_n^{-1}\Ext_{\mathcal{A}_*}(\F_p,H_*(Y))\]
so that the spectral sequence associated to $\tilde{v}^{-1}Y$ is the $v_n$-localised $\F_p$-ASS of $Y$.
\end{example}

\begin{example}
Let $X\in \Sp$. Set $Y_{-1}=\1\in\Syn_{\F_p}$ and assume by induction that we have defined $Y_{j-1}=:\1/(q_0^{i_0},\ldots,q_{j-1}^{i_{j-1}})\in\Syn_{\F_p}$ so that $\tau^{-1}Y_{j-1}=\Sph/(p^{i_0},\ldots,v_{j-1}^{i_{j-1}})$. Since everything is implicitly $p$-complete, $\1\in\Syn_{\F_p}$ is $\tau$-complete, hence so is $\mathrm{End}(Y_{j-1})$, and thus the corresponding $\sigma$-SS
\[
\pi_{*,*}(\mathrm{End}(Y_{j-1})/\tau)\implies \pi_*\mathrm{End}(\Sph/(p^{i_0},\ldots,v_{j-1}^{i_{j-1}}))
\]
converges. The self-map $v_j^{i_j}\in\pi_*\mathrm{End}(\Sph/(p^{i_0},\ldots,v_{j-1}^{i_{j-1}}))$ is therefore detected by an endomorphism $q_j^{i_j}$ of $Y_{j-1}/\tau$. Since $q_j^{i_j}$ is a permanent cycle in the $\sigma$-SS of $\mathrm{End}(Y_{j-1})$, it is also a permanent cycle in the $\tau$-BSS of $\mathrm{End}(Y_{j-1})$, and it therefore admits a lift to a self-map of $Y_{j-1}$ that $\tau$-inverts to $v_j^{i_j}$.

Denoting this self-map of $Y_{j-1}$ also by $q_j^{i_j}$, we set $Y_j=Y_{j-1}/q_j^{i_j}$. It follows that the synthetic spectrum 
\[\1/(q_0^{i_0},\ldots,q_n^{i_n}):=Y_n\]
has the property that 
\[\pi_{*,*}(\nu X\otimes\1/(q_0^{i_0},\ldots,q_n^{i_n},\tau))=\Ext_{\mathcal{A}_*}(\F_p,H(i_0,\ldots,i_n)\otimes H_*(X;\F_p))\]
so that the spectral sequence associated to $\nu X\otimes\1/(q_0^{i_0},\ldots,q_n^{i_n})$ is the modified $\F_p$-ASS of $X\otimes \1/(p^{i_0},\ldots,v_{n}^{i_{n}})$.
\end{example}


\section{Linear $t$-structures on synthetic spectra}\label{tstructuresection}

Let $E$ be a spectrum of Adams-type, such that $\Syn_E=\Syn$ is defined. In this section, we construct $t$-structures on $\Syn$ coming from lines in $\R^2$ through the origin. When restricted to the line defined by the equation $x=0$, we call the resulting $t$-structure the \emph{vertical $t$-structure}.


\subsection{Construction of the linear $t$-structure}

\begin{mydef}
    Let $L$ be a line in $\R^2$ through the origin not containing $(1,-1)$. Write $\R^2-L=L_{>}\cup L_{<}$ as the union of two connected subspaces where the region containing $(1,-1)$ is $L_>$. Let us also write $L_\geq=L\cup L_>$ and $L_\leq = L\cup L_<$. For any $n\in \Z$, define the $\infty$-subcategory $\Syn^L_{\geq n}$ of $\Syn$ to be spanned by those synthetic spectra $X$ such that $\pi_{a+n,b-n}X=\pi_{a,b}X[-n]=0$ if $(a,b)\in L_<$, so such that $X[-n]$ has homotopy groups supported in $L_{\geq}$. Likewise, define $\Syn^L_{\leq n}$ as the $\infty$-subcategory of $\Syn$ spanned by synthetic spectra $X$ with $\pi_{a-n,b+n}X=0$ for $(a,b)\in L_>$.\footnote{If $L$ is defined by $y=-x$, then there are two possible choices for $\Syn^{\leq L}$, and both could potentially produce $t$-structures on $\Syn$---we will not need either of them.}
\end{mydef}

\begin{theorem}\label{tstructuresexist}
    Let $L$ be a line in $\R^2$ through the origin not containing $(1,-1)$, $E$ be a spectrum such that $\Syn$ is cellular, and suppose that the unit $\1\in \Syn$ lies in $\Syn^L_{\geq 0}$. Then the $\infty$-subcategories $\Syn^L_\geq$ and $\Syn^L_\leq$ determine an accessible $t$-structure $\Syn$ which is both right and left complete and compatible with the monoidal structure on $\Syn$.
\end{theorem}

Note that Pstragowski's homological $t$-structure on $\Syn$ of \cite[\textsection4.2]{syntheticspectra} is not an instance of a linear $t$-structure, but there is a connection: by \cite[Cor.4.19]{syntheticspectra}, we see that a synthetic spectrum $X$ is connective in Pstragowski's $t$-structure if and only if $\nu E\otimes X$ is lies in $\Syn^L_{\geq0}$ for the line $L$ defined by $y=0$.\footnote{There are not many interesting $E$ where the line defined by $y=0$ produces a $t$-structure using \Cref{tstructuresexist}, as this would imply the $E$-ASS for $\Sph$ is concentrated in filtration $0$.}

\begin{proof}
    Consider the collection of objects $S=\{\Sigma^{a,b}\1\}$ where $(a,b)\in \L_\geq$. Define $\C_\geq$ as the full $\infty$-subcategory of $\Syn$ generated by the objects in $S$ by small colimits and extensions. By \cite[Pr.1.4.4.11(2)]{haname}, $\C_\geq$ is a presentable $\infty$-category, as $\Syn$ itself is presentable and stable \cite[Pr.4.2]{syntheticspectra}, and \cite[Pr.1.4.4.11(1)]{haname} then states that $\Syn$ has an accessible $t$-structure with with connective objects $\C_\geq$. Notice that by definition $\Syn^L_{\geq 0}$ contains the spheres in $S$. Moreover, it is clearly closed under extensions as this subcategory is defined by a vanishing condition on homotopy groups. Likewise, $\Syn^L_{\geq 0}$ is also closed under colimits, as coproducts and filtered colimits commute with homotopy groups and for pushouts we use the associated Mayer--Vietoris sequence on homotopy groups. It follows that $\C_{\geq}\subseteq \Syn_{\geq 0}^L$. The fact that $\Syn$ is cellular implies this inclusion is an equivalence of subcategories of $\Syn$.\\

To see this $t$-structure is left and right complete, if $\pi_{a,b}X=0$ for all $(a,b)$, then the synthetic spectrum $X=0$ from our cellularity hypothesis. In particular, we see that both of the $\infty$-subcategories $\bigcap \Syn^L_{\leq -n}$ and $\bigcap \Syn^L_{\geq n}$ contain only the zero synthetic spectrum. Just as in the proof of completeness in \cite[Pr.1.4.3.6]{haname}, referring to \cite[Pr.1.2.1.19]{haname} and its dual, it now suffices to show that $\Syn^L_{\leq 0}$ and $\Syn^L_{\geq 0}$ are closed under products and coproducts. This is clear though, as these categories are defined by a condition on synthetic homotopy groups, which commute with products, as well as finite coproducts and filtered colimits.\\
    
Finally, we want to check compatibility with the monoidal structure on $\Syn$. The unit lies in $\Syn_{\geq 0}^L$ by assumption, so it suffices to see that the objects of $S$ are closed under the tensor product of synthetic spectra. This is clear, as $\Sigma^{a,b}\1\otimes \Sigma^{c,d}\1\simeq \Sigma^{a+c,b+d}\1$, and if $(a,b),(c,d)\in L_\geq$, then also $(a+c,b+d)\in L_\geq$. 
\end{proof}

The condition about the unit $\1$ above depends on both $L$ and $E$. For example, when $E$ is connective and of Adams-type and $L$ is defined by $x=0$, which we call the \emph{vertical line}, then this unit condition is satisfied as the $E$-ASS for the $E$-local sphere are concentrated in nonnegative stems. For other $E$ and $L$, this relates to deep questions concerning vanishing lines in the $E$-ASS for the $E$-local sphere. The cellularity condition holds for all connective $\E_1$-rings of Adams-type by a recent result of Lawson; see \cite{lawson2024synthetic}. Our particular cases of interest are $E=\BP$ by a variant of \cite[Th.6.2]{syntheticspectra}, which proves the $E=\MU$-case, and for $E=\F_p$ by the following argument, which we include due to the importance of this example in this article.

\begin{lemma}\label{cellularityoffp}
The category $\Syn_{\F_p}$ is compactly generated by the bigraded spheres $S^{a,b}$ for any prime $p$.
\end{lemma}

\begin{proof}
By \cite[Rmk.4.14]{syntheticspectra}, $\Syn_{\F_p}$ is compactly generated by the set of bigraded suspensions of $\nu P$ where $P$ ranges through finite $\F_p$-projectives, so it remains to show $\nu P$ is in the subcategory $\mathcal{C}$ of $\Syn_{\F_p}$ generated by the bigraded spheres under colimits. Let $X$ be a bounded below spectrum, shifted so that its bottom homotopy group is in degree zero. We may form the following tower
\begin{equation}\label{fpresolutions}
\begin{tikzcd}
\bigoplus\limits_{\alpha_0}\Sph^0\arrow[d,"q_0"]&\bigoplus\limits_{\alpha_1}\Sph^1\arrow[d,"q_1"]&\bigoplus\limits_{\alpha_2}\Sph^2\arrow[d,"q_2"]&\\
X\arrow[r]&X_{(0)}\arrow[r]&X_{(1)}\arrow[r]&\cdots
\end{tikzcd}
\end{equation}
where $q_i$ is an isomorphism in $H_i(-;\Z)$, and $X_{(i)}$ is the cofibre of $q_i$, by use of the Hurewicz theorem. 

First, suppose $H_*(X;\F_p)=0$. It follows that $q_0$ induces the zero map in $H_*(-;\F_p)$ and $q_i$ induces an iso in $H_*(-;\F_p)$ for all $i>0$. Using \cite[Lm.4.23]{syntheticspectra}, the fact that $\colim X_{(i)}=0$, and that $\nu$ commutes with filtered colimits, one sees then that $\nu X\in\mathcal{C}$. If $P$ is any finite $\F_p$-projective, one has a localisation cofibre sequence
\[Q\to P\to P_{(p)}\]
Since $H_*(Q;\F_p)=0$, it follows again from \cite[Lm.4.23]{syntheticspectra} that
\[\nu Q\to \nu P\to\nu P_{(p)}\]
is a cofibre sequence in $\Syn_{\F_p}$. We may therefore assume that $P$ is a $p$-local finite spectrum.\\

Consider now a variation on (\ref{fpresolutions}), where we set $X=P$ and instead $q_i$ is chosen to be an isomorphism on $H_i(-;\F_p)$. The same argument now shows that $\nu X\in\mathcal{C}$ since $q_i$ induces a monomorphism on $H_*(-;\F_p)$ for all $i$.
\end{proof} 

\begin{remark}
Similar $t$-structures on filtered objects in a stable $\infty$-category were explored by Lee--Levy in \cite[\textsection2.2]{leelevi}. There is no moral difference between our construction and theirs---both are designed to cut out certain regions of spectral sequences. However, their construction takes place in filtered objects of a stable $\infty$-category, and ours takes place in synthetic spectra; we have not attempted to use the filtered model for $\Syn$ (\cite[Pr.C.22]{burkhahnsengartintate}) to relate our constructions. Moreover, we could consider $t$-structures on $\Syn$ defined with arbitrary pairs of subsets covering $\Z\times\Z$ satisfying certain axioms, similar to how Lee--Levy study the graphs of functions $\Z\to\Z$. For simplicity, we will ignore these generalisations here.
\end{remark}

\begin{mydef}
    Let $L$ and $E$ be as in \Cref{tstructuresexist} and consider the $t$-structure associated with this line on $\Syn$. Write $\tau^L_{\geq n}\colon \Syn\to \Syn^L_{\geq n}$ for the \emph{$n$th connective cover} functor and $\tau^L_{\leq n}\colon \Syn\to \Syn^L_{\leq n}$ for the \emph{$n$th truncation} functor. By \cite[1.4.4.13]{haname}, which states accessible $t$-structures lead to accessible connective cover and truncation functors, the adjoint functor theorem \cite[Cor.5.5.2.9]{httname}, and the presentability of $\Syn$ \cite[Pr.4.2]{syntheticspectra}, the functor $\tau^L_{\geq n}$ is right adjoint and $\tau^L_{\leq n}$ is left adjoint to their respective inclusions. If $L$ is the vertical line, then we write 
    \[ \Syn_{\geq n}^\uparrow=\Syn_{\geq n}^L, \qquad \Syn_{\leq n}^\uparrow=\Syn_{\geq n}^L, \qquad \tau_{\geq n}^L=\tau_{\geq n}^\uparrow, \qquad \tau_{\leq n}^L=\tau_{\leq n}^\uparrow\]
    and refer to the associated $t$-structure as the \emph{vertical $t$-structure}.
\end{mydef}


\subsection{The heart of the linear $t$-structure}

To calculate the heart of these $t$-structures, we need the following algebraic construction.

\begin{mydef}
    Let $L$ and $E$ be as in \Cref{tstructuresexist}. For $X\in\Syn$, write $\pi^L_\ast X$ for the bigraded homotopy groups $\pi_{a,b} X$ where $(a,b)\in L$. As $L$ is a line through the origin and the $t$-structure defined by $L$ is monoidal, then $\pi^L_\ast\1$ is a graded ring.
\end{mydef}

\begin{prop}\label{calculationoftstructures}
    Let $L$ and $E$ be as in \Cref{tstructuresexist}. Then the functor
    \[\pi_\ast^L\colon (\Syn^L)^\heartsuit \xrightarrow{\simeq} \Mod^\gr_{\pi_\ast^L \1}\]
    is an equivalence of $\infty$-categories between the heart of the $t$-structure associated with $L$ and the abelian category of graded modules over $\pi_\ast^L \1$.
\end{prop}

\begin{proof}
The functor 
\[\pi^L_\ast\colon (\Syn^L)^\heartsuit\to \Mod^\gr_{\pi_\ast^L \1}\]
is conservative as $\Syn_{\geq 0}^L$ is generated by the bigraded spheres $\Sigma^{a,b}\1$ with $(a,b)\in L_{\ge}$ by our cellularity assumption. Moreover, it preserves cokernels, as for any $(a,b)\in L$ and any $Y\in\Syn$,
\[[\Sigma^{a,b}\1,Y]=[\Sigma^{a,b}\1,\tau_{\le0}^L Y]\]
and if $X\in \Syn_{\geq 0}^L$, then
\[[\Sigma^{a,b}\1,\Sigma X]=[\Sigma^{a-1,b+1}\1,X]=0\]
as $(a-1,b+1)\in L_{<}$. A similar argument shows that $\pi^L_\ast$ preserves kernels and products, and it follows that $\pi^L_\ast$ has a left adjoint $H(-)$.\\

The classical Barr--Beck theorem implies that the adjunction is monadic, and it remains to show the monad $\pi_*^LH(-)$ is the identity functor. Using that $\pi_*^L$ commutes with filtered colimits, finite coproducts, and cokernels, it suffices to show that the counit map $\epsilon\colon H(\pi_*^L\1)\to\tau_{\le0}^L\1$ is an equivalence. For any $X\in (\Syn^L)^\heartsuit$, this induces a map
\[[\tau_{\le0}^L\1,X]\to[H(\pi_*^L\1),X]\]
Both sides are isomorphic to $\pi_{0,0}X$, on the left because $X\simeq \tau_{\le0}^LX$ and on the right by adjunction. Using generation by bigraded spheres, it suffices to check this map is an isomorphism when $X=\tau_{\le0}^L\Sigma^{a,b}\1$ for all $(a,b)\in L$, which together give a map of $\pi_*^L\1$-modules $\pi_*^L\1\to\pi_*^L\1$, which is an isomorphism because it sends 1 to the map adjoint to the counit $\epsilon$, in other words, the identity. 
\end{proof}

\begin{remark}
Although two different lines $L$ and $L'$ produce different $t$-structures on $\Syn$, it follows from \Cref{calculationoftstructures}, that many $t$-structures give the same heart. For example, all lines $L$ defined by the equation $y=ax$, where $a$ is irrational, have heart equivalent to graded modules over $\pi_\ast^L\1\simeq \pi_{0,0}\1$.
\end{remark}


\subsection{Homotopy groups of linear connective covers and truncations}

Let us end this section by calculating some $L$-connective covers and $L$-truncations, as well as the compatibility of these $t$-structures with the usual $t$-structure on spectra.

\begin{prop}\label{htpyofcovers}
    Let $L$ and $E$ be as in \Cref{tstructuresexist} and $X\in\Syn$. Furthermore, suppose that $L$ has slope $m$ with $-1<m\le \infty$. Then we have natural isomorphisms of abelian groups
    \[\pi_{a,b} \tau_{\geq n}^L X\simeq \begin{cases}
        \pi_{a,b} X & (a-n,b+n)\in L_\geq   \\
        0 & (a-n,b+n)\in L_<
    \end{cases}\]
    For the line defined by $x=0$, we have natural isomorphisms of abelian groups
    \[\pi_{a,b} (\tau_{\geq n}^\uparrow(X)/\tau)\simeq \begin{cases}
        \pi_{a,b}  (X/\tau) & a> n   \\
        Z^\infty & a=n    \\
        0 & a< n.
    \end{cases}\]
    where $Z^\infty$ is the subgroup of $\pi_*^\uparrow(X/\tau)$ consisting of permanent cycles in the $\tau$-BSS for $X$.
\end{prop}

\begin{remark}
    Computing $\tau_{\geq n}^L X/\tau$ for general $L$ is subtle. What simplifies the expression in the case of the vertical $t$-structure is the fact that there cannot be differentials which cross this line for degree reasons; only with source or target on the line. In the case of a more general line $L$, the appearance of $Z_\infty$ above needs to be appropriately replaced by kernels of \emph{truncated total differentials}; see \cite[Df.2.19]{smfcomputation}. 
\end{remark}

\begin{proof}[Proof of \Cref{htpyofcovers}]
The first claim follows from the definitions of $\Syn^L_{\ge n}$ and $\Syn^L_{\le n-1}$ and the long exact sequence associated with
\begin{equation}\label{Lcoverfibre}\tau_{\ge n}^LX\to X\to \tau_{\le n-1}^LX\end{equation}
The claim about homotopy groups mod $\tau$ follows from investigating the diagram of long exact sequences
\begin{equation}\label{taules}
\begin{tikzcd}
\pi_{a,b+1}\tau_{\ge n}^LX\arrow[r,"\tau"]\arrow[d]&\pi_{a,b}\tau_{\ge n}^LX\arrow[r]\arrow[d]&\pi_{a,b}(\tau_{\ge n}^L(X)/\tau)\arrow[r, "\partial"]\arrow[d]&\pi_{a-1,b+2}\tau_{\ge n}^LX\arrow[r,"\tau"]\arrow[d]&\pi_{a-1,b+1}\tau_{\ge n}^LX\arrow[d]\\
\pi_{a,b+1}X\arrow[r,"\tau"]&\pi_{a,b}X\arrow[r]&\pi_{a,b}(X/\tau)\arrow[r,"\partial"]&\pi_{a-1,b+2}X\arrow[r,"\tau"]&\pi_{a-1,b+1}X
\end{tikzcd}
\end{equation}
If $(a-n,b+n)\in L$, then the second vertical map from the left is an isomorphism, and the last two terms in the top row are zero. If $L$ is not vertical, then the top left term is also zero, and the claim follows in this case. If $L$ is vertical, the leftmost vertical map is an isomorphism, and this identifies $\pi_{a,b}(\tau_{\ge n}^L(X)/\tau)$ with the kernel of $\partial$, which is the subgroup of permanent cycles in the $\tau$-BSS.
\end{proof}

\begin{remark}\label{truncationsswitcharoo}
Similar claims hold for the $L$-truncations using (\ref{Lcoverfibre}) and (\ref{taules}).
\end{remark}

\begin{remark}
One can make a more general statement when the line $L$ has slope $m$ with $-\infty<m<-1$. For example, the proof of \Cref{htpyofcovers} shows that if $m=-r$ for some integer $2\leq r<\infty$, then the $Z_L^\infty$ term in the case of $(a-n,b+n)\in L$, needs to be replaced with $Z_L^{r-1}$, the subgroup of $\pi_\ast^L(X/\tau)$ consisting of $d_{\le r-1}$-cycles in the $\tau$-BSS for $X$.
\end{remark}

\begin{cor}\label{coverstaucomplete}
Let $L$ and $E$ be as in \Cref{htpyofcovers}, $X$ be a $\tau$-complete synthetic spectrum, and $n$ be an integer. Then both $\tau_{\ge n}^L X$ and $\tau_{\le n}^L X$ are also $\tau$-complete.
\end{cor}

\begin{proof}
We will prove that $\tau_{\ge n}^L X$ is $\tau$-complete for all $n$; the same follows for $\tau_{\le n}^L X$ since being $\tau$-complete is closed under taking cofibres. Consider the diagram of Milnor sequences
\begin{equation}
\begin{tikzcd}
0\arrow[r]&\lim^1_m\pi_{a+1,b-1}(X_n/\tau^m)\arrow[r]\arrow[d]&\pi_{a,b}(\lim_m X_n/\tau^m)\arrow[r]\arrow[d]&\lim_m\pi_{a,b}(X_n/\tau^m)\arrow[r]\arrow[d]&0\\
0\arrow[r]&\lim^1_m\pi_{a+1,b-1}(X/\tau^m)\arrow[r]&\pi_{a,b}X\arrow[r]&\lim_m\pi_{a,b}(X/\tau^m)\arrow[r]&0
\end{tikzcd}
\end{equation}
where $X_n=\tau_{\ge n}^L(X)$. Applying \Cref{htpyofcovers}, we see that if $(a-n,b+n)\in L_{>}$, then there is an isomorphism $\pi_{a,b}\tau_{\ge n}^L(X)\simeq \pi_{a,b} X$ and $\pi_{a+1,b-1}(\tau_{\ge n}^L(X))\simeq \pi_{a+1,b-1} X$. The diagram above then shows that
\[\pi_{a,b}\tau_{\ge n}^L(X)\to\pi_{a,b}(\lim_m\tau_{\ge n}^L(X)/\tau^m)\]
is an isomorphism. If $(a-n,b+n)\in L_{<}$, then 
\[\pi_{a,b}(\lim_m\tau_{\ge n}^L(X)/\tau^m)=\lim_m\pi_{a,b}(\tau_{\ge n}^L(X)/\tau^m)=\pi_{a,b}\tau_{\ge n}^L(X)=0\]
so it remains to show that $\lim^1_m\pi_{a+1,b-1}(\tau_{\ge n}^L(X)/\tau^m)=0$. This is automatic if $(a+1-n,b-1+n)\in L_{<}$, and if $(a+1-n,b-1+n)\in L_{\ge}$, this follows as one has $\pi_{s,t}(\tau_{\ge n}^L(X)/\tau^m)=\pi_{s,t}(\tau_{\ge n}^L(X))/\tau^m$ whenever $(s-n,t+n)\in L_{\ge}$ and $(s-1-n,t+1+n)\in L_{<}$, so that the inverse system is Mittag-Leffler.\\

If $(a-n,b+n)\in L$, the lefthand vertical map is an isomorphism as before, and the result follows from the claim that the map
\[\lim_m\pi_{a,b}(\tau_{\ge n}^L(X)/\tau^m)\to\lim_m\pi_{a,b}(X/\tau^m)\]
is an isomorphism. To see this, first note that for each $m$
\[\pi_{a,b}(\tau_{\ge n}^L(X)/\tau^m)\to\pi_{a,b}(X/\tau^m)\]
is an injection as it is identified with the map $\pi_{a,b}(X)/\tau^m\to \pi_{a,b}(X/\tau^m)$, which implies the map given by taken limits is also an injection. The map also becomes a surjection in the limit because, by $\tau$-completeness of $X$, each class on the right-hand side admits a lift to $\pi_{a,b}X$, whose projection to $\pi_{a,b}(X/\tau^m)$ must factor through $\pi_{a,b}(X)/\tau^m$ since the composite $\pi_{a,b}X\to \pi_{a,b}(X/\tau^m)\xrightarrow{\partial}\pi_{a-1,b+1}X$ is zero.
\end{proof}

Although we have worked mostly with a general line through the origin in this section, there are some statements that behave best with respect to vertical $t$-structure. For example, the following interaction between vertical connective covers and truncations and $\tau$-inversion is an immediate consequence of \Cref{htpyofcovers}.

\begin{cor}\label{tauinversionofcovers}
   Let $L$ be the vertical line, $E$ be as in \Cref{tstructuresexist}, $X\in\Syn$, and $n$ be an integer. Then there are natural equivalences of spectra
    \[\tau^{-1}\tau_{\geq n}^\uparrow X\simeq \tau_{\geq n} \tau^{-1}X\qquad\qquad \tau^{-1}\tau_{\leq n}^\uparrow X\simeq \tau_{\leq n} \tau^{-1}X.\]
    In particular, if $Y\in \Sp$ is connective, then $\nu Y$ is vertically connective.
\end{cor}

\begin{proof}
Since $\tau$ has bidegree $(0,-1)$, we have that
\[\tau^{-1}\pi_{k,*}\tau_{\ge n}^\uparrow X=\pi_k(\tau^{-1}X)[\tau^{\pm}]\]
for all $k$, so that in particular the homotopy groups of $\tau^{-1}\tau_{\geq n}^\uparrow X$ are concentrated in degrees $\geq n$. In particular, the natural map $\tau^{-1}\tau_{\geq n}^L X\to \tau^{-1} X$ factors through a map $\tau^{-1}\tau_{\geq n}^L X\to \tau_{\geq n} \tau^{-1}X$, which is an equivalence by \Cref{htpyofcovers}. For truncations, combine the same argument with \Cref{truncationsswitcharoo}.\\

The ``in particular'' statement follows this together with \Cref{filspthm} applied to $\sigma \nu Y$.
\end{proof}

\begin{remark}
To see the necessity of $L$ being the vertical line, consider $E=\BP$ at the prime $2$, $X=\j_\BP$ as defined in \Cref{imageofjsyntheticdefinition} to come, and let $L$ be defined by the equation $x=y$. Then the connective cover of $\tau^{-1}\j_\BP$ is equivalent to $\j$, and in particular, $\pi_4\tau^{-1}\j_{\BP}=0$. On the other hand, from \Cref{manssjattwopictureetwo}, we see that $\pi_4 \tau^{-1}\tau_{\geq0}^L \j_\BP\neq 0$. Indeed, the connective cover will cut off exactly the top-diagonal line of red circles in the $\sigma$-SS associated with $\tau_{\geq 0}^L\j_\BP$, meaning the class in degree $(4,2)$ on the $E_2$-page cannot support any differentials, so it survives to a nonzero class in $\pi_4 (\tau^{-1}\tau_{\geq0}^L \j_\BP)$. In general, if $X$ is an $E$-nilpotent complete spectrum with the property that the $E_2$-page of its $E$-ASS vanishes in $L_<$, then $\tau^{-1}\tau_{\geq 0}^L \nu X\simeq \tau_{\geq0} \tau^{-1} \nu X$. Similar statements hold for certain vanishing \emph{curves} with associated $t$-structure given by the graph of said curve; we will come back to this again in future work.
\end{remark}


\section{Synthetic analogues of topological $K$-theory}\label{ktheoryconstructionsection}

As a warm-up to our synthetic constructions and calculations of various image-of-$J$ spectra, we begin with an analysis of topological $K$-theories in synthetic spectra; recall that we are only starting with the basic assumptions of \Cref{sec:basicassumptions}.\\

From the point of view of this article, it is not necessary to always consider synthetic analogues such as $\nu \ko$. Often, we make a purely synthetic construction, such as that of $\KO_{\BP}$ in \Cref{woodstheoremsection}, show it is a synthetic lift of $\KO$ and compute the associated spectral sequence, and then remark that it actually computes the ANSS for $\KO$. In other words, we do not need to label a modified $E$-ASS as a classical $E$-ASS to work with it; a purely synthetic perspective suffices as the end of the day to compute spectral sequences for $\ko$ and $\j$ to prove \Cref{maintheorem}.\\

We begin with a synthetic proof of the classical \emph{Wood's theorem}. Let us remind the reader that we are implicitly working in the $p$-complete setting.


\subsection{Wood's theorem}\label{woodstheoremsection}
We already know the homotopy groups of $\KU$ and we have defined $\KO=\KU^{hC_2}$. Our next goal is to show how synthetic spectra can be used to compute the homotopy groups of $\KO$ (\Cref{gradedringofktheory}) as well as prove the following classical results. For this subsection, we set $p=2$.

\begin{theorem}[Wood's theorem {\cite{wood}}]\label{theclassicalwoodstheorem}
    Let $\eta \in \pi_1 \KO$ denote the image of the Hopf element $\eta \in \pi_1 \Sph$. Consider the composition
    \begin{equation}\label{woodcofibresequence}\Sigma \KO\xrightarrow{\eta} \KO\xrightarrow{c} \KU.\end{equation}
    where $c$ is the canonical projection. Then (\ref{woodcofibresequence}), which represents the zero map as it corresponds to a class in $\pi_1 \KU$ as a map of $\KO$-modules, is a cofibre sequence. We call this the \emph{Wood cofibre sequence}.
\end{theorem}

Along to way to prove this theorem, we will also compute the homotopy groups of $\KO$, which leads to the connective version of Wood's theorem.

\begin{cor}\label{cor:connectivewoodstheorem}
    Let $\eta \in \pi_1 \ko$ denote the image of the Hopf element $\eta \in \pi_1 \Sph$. Then the composition
    \[\Sigma \ko\xrightarrow{\eta} \ko\xrightarrow{c} \ku\]
    is a cofibre sequence of $\ko$-modules.
\end{cor}

\begin{proof}
    The statement of \Cref{theclassicalwoodstheorem} is equivalent to an equivalence of $\KO$-modules $\KO\otimes C(\eta) \simeq \KU$. Taking connective covers we get an equivalence of $\ko$-modules $\tau_{\geq 0} (\KO\otimes C(\eta)) \simeq \ku$. To prove \Cref{cor:connectivewoodstheorem}, it then suffices to show that $\ko \otimes C(\eta)$ is equivalent to $\tau_{\geq 0} (\KO\otimes C(\eta))$, which follows from the fact that $\pi_d\KO=0$ for $-3\leq d\leq -1$; see \Cref{gradedringofktheory} below.
\end{proof}

There are many proofs of Wood's theorem, including the classical reference \cite{wood}. If one uses our set of assumptions from \Cref{sec:basicassumptions}, then Mathew's approach in the first arXiv version of \cite[Th.3.2]{tmfhomology} is a slick proof entirely within the realm of chromatic homotopy theory. Our proof of \Cref{theclassicalwoodstheorem} acts as a prelude to many of the techniques used in the rest of this article and bears a similarity to an approach using the ANSS.\\

Let us begin by defining a $\BP$-synthetic lift of $\KO$.

\begin{mydef}\label{definitionofko}
    Let $\KO_\BP = (\nu\KU)^{hC_2}$ in $\Syn_{\BP}$.
\end{mydef}

\begin{lemma}
    There is an equivalence of $\E_\infty$-rings $\tau^{-1} \KO_\BP \simeq \KO$.
\end{lemma}

\begin{proof}
This is \cite[Pr.1.6]{osyn}, and the argument is simple: first, recall the homological $t$-structure on $\Syn_{\BP}$ from \cite{syntheticspectra}, and let $\tau_{\ge0}^{\mathrm{hom}}$ denote the corresponding connective cover and $(\Syn_{\BP})_{\ge0}$ the connective objects. The functor $\nu:\Sp\to(\Syn_\BP)_{\ge0}$ is right adjoint to the restriction of $\tau^{-1}$ to $(\Syn_\BP)_{\ge0}$ by \cite[Th.4.37 \& Pr.4.33]{syntheticspectra}, and $\tau^{-1}$ vanishes on coconnective objects in $\Syn_\BP$ by \cite[Lm.4.35]{syntheticspectra}. It follows that
\[\tau^{-1}(\nu(\KU)^{hC_2})=\tau^{-1}(\tau_{\ge0}^{\mathrm{hom}}(\nu(\KU)^{hC_2}))=\tau^{-1}(\tau_{\ge0}^{\mathrm{hom}}(\nu(\KU^{hC_2})))=\tau^{-1}(\nu\KO)=\KO.\qedhere\]
\end{proof}

\begin{remark}
    More is true. In fact, $\KO_\BP \simeq \nu \KO$, so we should not be surprised that the spectral sequence associated to $\KO_\BP$ looks like the ANSS for $\KO$. However, we will not use this fact. We will prove this in \Cref{agreementinbpcase} by direct computation, but it is also true for purely abstract reasons by \cite[(0.1) \& Th.A]{osyn}.
\end{remark}

To prove \Cref{theclassicalwoodstheorem}, we will need a few synthetic lemmata. The key step is the detection lemma; see \Cref{etaisdetectionmodtau}.

\begin{lemma}\label{etwoforko}
    The spectral sequence associatied to $(\nu \KU)^{hC_2}$ has $E_2$-page of the form
    \[E_2^{s,f} = \pi_{s,f} (\nu \KU)^{hC_2}/\tau \simeq \Z[u^2, x]/(2x) \implies \pi_s \tau^{-1}((\nu \KU)^{hC_2}) = \pi_s \KO,\]
    where $|x| = (1,1)$ and $|u^2| = (4,0)$.
\end{lemma}

\begin{proof}
    The homotopy fixed point spectral sequence (HFPSS) internal to $\Syn_{\BP}$ applied to $(\nu \KU)^{hC_2}/\tau \simeq (\nu \KU/\tau)^{hC_2}$ takes the form
    \[E_2 \colon H^t(C_2; \pi_{s,f} \nu\KU/\tau) \implies \pi_{s-t,f+s} (\nu \KU/\tau)^{hC_2},\]
    where differentials $d_r$ shift the degree by $(-1,+1,+r)$. It is easy to compute the effect of the $C_2$-action on $\pi_{\ast,\ast}\nu\KU/\tau$ by \Cref{examplesofsignatures}, as we know how $C_2$-acts on the $E_2$-page of the ANSS for $\KU$. In particular, $C_2$ acts trivially on $\pi_{s,\ast}\nu\KU$ when $s$ is divisible by $4$, and by the sign action when $s\equiv 2$ modulo $4$. A standard group cohomology calculation then gives the claimed $E_2$-page.
    There are no nonzero differentials in this spectral sequence for degree reasons.
\end{proof}

\begin{lemma}\label{etaisdetectionmodtau}
    Let $E=\BP$ and $p=2$. The element $h_1 \in \pi_{1,1} \1/\tau$ which detects $\eta\in \pi_1 \Sph$ has nonzero image $x$ in $\pi_{1,1} (\nu \KU)^{h C_2}/\tau$.
\end{lemma}

\begin{proof}
    First, note that the cofibre sequence
    \[\Sph \xrightarrow{2} \Sph \to \Sph/2\]
    induces a short exact sequence on $\BP_\ast$-homology, so $\nu(\Sph/2) \simeq \1/2$; see \cite[Lm.4.23]{syntheticspectra}. We can then identify the spectral sequence associated to $\1/2$ with the ANSS of $\Sph/2$ by \Cref{examplesofsignatures}. We then have the computation
    \[\pi_{\ast,0} \1/(\tau,2) \simeq \Ext_{\BP_*\BP}^{0,*}(\BP_*,\BP_*/2)=\F_2[v_1];\]
    see \cite[Th.4.3.2]{greenbook}. Since $2\cdot h_1=0\in\pi_{1,1}\1/\tau$, exactness implies that $h_1=\partial(v_1)$, where $\partial \colon \1/(\tau,2) \to \Sigma^{1,-1}\1/\tau$ is the boundary map from the defining cofibre sequence of $\1/(\tau,2)$.\\

    As $\KU$ comes equipped with the multiplicative formal group law, we know that $v_1 \in \pi_2 \BP/2$ is mapped to the nonzero element $\overline{u}$ in $\pi_2 \KU/2$. This means that right unit $\eta_R\colon \BP_\ast/2 \to \BP_\ast \KU/2$ sends $v_1$ to $\overline{u}$. On the other hand, the left unit $\eta_L$ sends $v_1$ to $v_1$ by definition. By \cite[(4.3.1)]{greenbook}, one has $\eta_L(v_1)\equiv\eta_R(v_1)$ modulo $2$. We conclude that the image of $v_1$ under the unit $\1/(\tau,2) \to \nu\KU/(\tau,2)$ is nonzero.\\
    
    In particular, this means that the image of $v_1$ in $\pi_{2,0} (\nu \KU)^{h C_2}/(\tau,2)$ is also nonzero, as the unit map $\1 \to \nu\KU$ factors through $(\nu \KU)^{h C_2}$. Consider the commutative diagram of abelian groups
    \[\begin{tikzcd}
            &   {\pi_{2,0}\1/(\tau,2)}\ar[d, "{\overline{h}}"]\ar[r, "\partial"]    &   {\pi_{1,1} \1/\tau}\ar[d, "{h}"]    \\
        {\pi_{2,0} (\nu \KU)^{h C_2}/\tau}\ar[r]    &   {\pi_{2,0}(\nu \KU)^{h C_2}/(\tau,2)}\ar[r, "\partial"]    &   {\pi_{1,1} (\nu \KU)^{h C_2}/\tau}
    \end{tikzcd}\]
    where the rows are exact and the vertical maps are unit maps.  Since $\overline{h}(v_1)$ is nonzero and $\partial(v_1)=h_1$, to see that $h(h_1)\neq0$, it suffices to show that the lower $\partial$ map is injective. This is because $\pi_{2,0}(\nu\KU)^{hC_2}/\tau=0$ by \Cref{etwoforko}.
\end{proof}

\begin{lemma}\label{signatureofourstupidko}
    The spectral sequence associated to $(\nu \KU)^{hC_2}$ is determined as a multiplicative spectral sequence by the differential $d_3(u^2) = h_1^3$.
\end{lemma}

See \Cref{mansskoattwopicture} for a connective version of this spectral sequence.

\begin{proof}
    The unit map $\1 \to (\nu\KU)^{hC_2}$ induces a map from the ANSS for $\Sph$ to the spectral sequence associated to $(\nu\KU)^{hC_2}$. In particular, as we have assumed knowledge of the ANSS for $\Sph$ in degrees $\leq 8$, we see that the class $h_1 \in \pi_{1,1}\1/\tau$ has the property that $h_1^4$ is hit by a $d_3$. This means that $\tau^2h_1^4=0 \in \pi_{4,2}\1$ by \Cref{filspthm}, and hence $\tau^2h_1^4=0 \in \pi_{4,2}(\nu\KU)^{hC_2}$. By \Cref{etaisdetectionmodtau}, the image of $h_1^4$ is nonzero in $\pi_{4,4}(\nu\KU)^{hC_2}/\tau$, so $h_1^4$ must be hit by a $d_{
    \le 3}$. The Leibniz rule then implies the desired differential $d_3(u^2) = h_1^3$. The Leibniz rule then propagates and the spectral sequence then collapses on the $E_4$-page.
\end{proof}

\begin{lemma}\label{gradedringofktheory}
    There is an isomorphism of graded rings
    \[\pi_\ast \KO \simeq \Z[\eta, \al, \be^\pm]/(2\eta, \eta^3, \eta\al, 4\be-\al^2),\qquad |\eta|=1, |\al|=4, |\be|=8.\]
    The projection map $c\colon \KO \to \KU$, induced from the equivalence $\tau^{-1} \nu \KU \simeq \KU$, can be described on homotopy groups by the formul{\ae}
    \[c(\eta)=\eta, \qquad c(\al)= 2u^2, \qquad c(\be) = u^4,\]
    where $u\in \pi_2 \KU$ is the generator.
\end{lemma}

\begin{proof}
    This follows immediately from \Cref{signatureofourstupidko} and the associated map of spectral sequences from the spectral sequences associated to $(\nu \KU)^{hC_2}$ to the ANSS of $\KU$.
\end{proof}

\begin{lemma}\label{sigatureofmodtwoktheory}
    The spectral sequence associated to $(\nu\KU/2)^{hC_2} \simeq (\nu\KU)^{hC_2}/2$ takes the form
    \[E_2^{s,f} = \pi_{s,f}(\nu\KU/2)^{hC_2}/\tau \simeq \F_2[\overline{u}^\pm, h_1] \implies \pi_s \KO/2, \]
    where $|h_1|=(1,1)$, $|\overline{u}|=(2,0)$. It is determined as a module over the spectral sequence associated to $(\nu\KU)^{hC_2}$ by the fact that $1$ and $\overline{u}$ are permanent cycles.
\end{lemma}

We write the $E_2$-page as a ring for simplicity, even though this spectral sequence has no multiplicative structure.

\begin{proof}
    First, notice that $\tau$-inversion commutes with cofibres, so the above spectral sequence has the correct abutment. The $E_2$-page is computed as in the proof of \Cref{etwoforko}, although now the $C_2$-action on $\pi_{\ast,\ast}\nu\KU/(2,\tau)$ is trivial. The $E_2$-page is then generated as a module over the $E_2$-page of the spectral sequence associated to $(\nu\KU)^{hC_2}$ by the classes $1$ and $\overline{u}$. Since $1$ is a permanent cycle in the spectral sequence associated to $(\nu\KU)^{hC_2}$, the fact that the map of spectral sequences induced by the map $(\nu\KU)^{hC_2}\to(\nu\KU/2)^{hC_2}$ sends $1$ to $1$ implies the class $1$ is a permanent cycle in the target. The proof of \Cref{etaisdetectionmodtau} shows that the map $\pi_{2,0}\1/(\tau,2)\to\pi_{2,0}(\nu\KU/2)^{hC_2}/\tau$ sends $v_1\mapsto \overline{u}$. Since $v_1$ is a permanent cycle in the source, $\overline{u}$ is thus a permanent cycle in the target.\\
    
    The module structure over the spectral sequence associated to $(\nu\KU)^{hC_2}$ determines all $d_3$-differentials, after which the spectral sequence collapses for degree reasons.
\end{proof}

\begin{lemma}\label{homotopyofmodtwoktheory}
    The homotopy groups of $\KO/2$ as given by
    \[\pi_{s} \KO/2 \simeq \begin{cases}
        \F_2\{\overline{u}^{s/2}\}   &   {s\equiv_8 0}   \\
        \F_2\{\eta\overline{u}^{(s-1)/2}\}   &   {s\equiv_8 1}   \\
        \Z/4\Z\{\overline{u}^{s/2}\}   &   {s\equiv_8 2}   \\
        \F_2\{\eta\overline{u}^{(s-1)/2}\}   &   {s\equiv_8 3}   \\
        \F_2\{\eta^2\overline{u}^{(s-2)/2}\}   &   {s\equiv_8 4}   \\
        0   &   {s\equiv_8 5,6,7,}
    \end{cases}\]
    with the hidden extension $\eta^2\overline{u}^{4k} = 2\overline{u}^{4k+1}$.
\end{lemma}

\begin{proof}
    This follows straight from \Cref{sigatureofmodtwoktheory} as the associated spectral sequence has a horizontal vanishing line. The interesting extension in degrees congruent to $2$ modulo $8$ comes from the classical extension $2 v_1 = \eta^2$ from the sphere, which is one of our basic assumptions of \Cref{sec:basicassumptions}.
\end{proof}

\begin{proof}[Proof of \Cref{theclassicalwoodstheorem}]
    The goal is to identify the cofibre of $\eta \colon \Sigma \KO \to \KO$ with $\KU$ in the $\infty$-category of $\KO$-modules. From this we attempt to compute the homotopy groups of $C(\eta)\otimes \KO$ using the associated long exact sequence on homotopy groups. This is simple, except for $\pi_{8d+4} C(\eta)\otimes \KO$, where we have to solve the extension problem
    \[0\to \Z \to \pi_{8d+4} C(\eta)\otimes \KO \to \F_2 \to 0.\]
    By periodicity, it suffices to prove the $d=0$ case. If this extension was split, then $\pi_4 (C(\eta)\otimes \KO)/2$ would be a $2$-dimensional $\F_2$-vector space. However, we have an exact sequence 
\[\pi_3\KO/2\xrightarrow{\eta, \simeq}\pi_4\KO/2\to\pi_4(C(\eta)\otimes\KO)/2\to \pi_2\KO/2\xrightarrow{\eta}\pi_3\KO/2, \]    
    
    and using \Cref{homotopyofmodtwoktheory}, we conclude that $\pi_4(C(\eta)\otimes \KO)/2\simeq \F_2$, hence the above extension is not split. In particular, we see that the homotopy groups of $C(\eta)\otimes \KO$ are isomorphic to the homotopy groups of $\KU$. We have a natural map $C(\eta)\otimes \KO \to \KU$ from the fact that the composite (\ref{woodcofibresequence}) is the zero map. Comparing this map with the computations of the map $c$ in \Cref{gradedringofktheory} and the computations of $C(\eta)\otimes \KO$ above, we obtain the desired conclusion.
\end{proof}

For later use, we include a corollary of our computations above.

\begin{cor}\label{cor:woodboundary}
    Writing $\partial\colon \KU \to \Sigma^2\KO$ for the boundary map of the cofibre sequence (\ref{woodcofibresequence}), then for $k\in \Z$ we have
    \[\partial(u^{4k}) = 0 , \qquad \partial(u^{4k+1}) = \pm2\be^k , \qquad \partial(u^{4k+2}) = \eta^2\be^k , \qquad \partial(u^{4k+3}) = \pm\al\be^k,\]
    using the notation of \Cref{gradedringofktheory}.
\end{cor}

\begin{proof}
    This follows immediately from the $\eta$-action on $\pi_\ast \KO$ and the formul{\ae} for $c$, both of which appear in \Cref{gradedringofktheory}.
\end{proof}


\subsection{Synthetic Wood's theorem}\label{syntheticwoodstheoremsection}

Having proven the classical Wood's theorem in the previous subsection, we can now upgrade this to a synthetic statement.

\begin{prop}\label{synwood}
There is a cofibre sequence in $\Syn_{E}$
\[\Sigma^{1,1}\nu\ko\xrightarrow{h_1}\nu\ko\xrightarrow{\nu c}\nu\ku\]
for $E=\F_2$ and $E=\BP$, which is a synthetic lift of the cofibre sequence of \Cref{cor:connectivewoodstheorem}, where $h_1$ is the unique class in $\pi_{1,1}\1$ such that $\tau^{-1}(h_1)=\eta$.
\end{prop}

\begin{proof}
The cofibre sequence of spectra
\[\Sph^0\to C(\eta)\to\Sph^2\]
induces a short exact sequence on $\F_2$- and $\BP$-homology, so by \cite[Lm.4.23]{syntheticspectra}, it follows that one has a cofibre sequence
\[\1\to \nu C(\eta)\to \Sigma^{2,0}\1\]
in $\Syn_E$. The $\tau$-inversion is an exact functor, so we see that the boundary map $\partial\colon \Sigma^{1,1}\1\to \1$ in this cofibre sequence must $\tau$-invert to $\eta$. This $\partial$ then represents the desired $h_1\in \pi_{1,1}\1$. In particular, we have an equivalence of synthetic spectra $\nu C(\eta)\simeq C(h_1)$. As $C(\eta)$ is both an $\F_2$- and $\BP$-finite projective, we use \cite[Lm.4.24]{syntheticspectra} to see that
\[\nu\ku\simeq \nu(\ko\otimes C(\eta))\simeq \nu\ko\otimes \nu C(\eta)\simeq \nu\ko\otimes C(h_1)\]
by \Cref{cor:connectivewoodstheorem}.
\end{proof}

\begin{warn}
According to the last line in the above proof, there is an equivalence of $\F_2$-synthetic spectra $\nu\ku\simeq \nu\ko\otimes C(h_1)$. Not all classical relationships between $\ku$ and $\ko$ lift to synthetic spectra though. For example, $\tau_{\geq 0}^\uparrow (\nu \ku)^{hC_2}$ is \textbf{not} equivalent to $\nu\ko$; see \Cref{koF2nogo}.
\end{warn}


\subsection{$\nu\ku$ and $\nu\ell$ in $\F_p$-synthetic spectra}
We start with $\ku$ and $\ell$.

\begin{prop}\label{synthetichomotopyofkuandell}
In $\Syn_{\F_2}$, we have
\[\pi_{*,*}\nu \ku=\Z[h_0,q_1,\tau]/(\tau h_0-2)\]
where $|h_0|=(0,1)$, $|q_1|=(2,1)$, and $|\tau|=(0,-1)$. In $\Syn_{\F_p}$ for $p>2$, we have 
\[\pi_{*,*}\nu \ell=\Z[h_0,q_1,\tau]/(\tau h_0-p)\]
where $|h_0|=(0,1)$, $|q_1|=(2p-2,1)$, and $|\tau|=(0,-1)$.
\end{prop}

\begin{proof}
    Let us discuss the $p=2$ case for $\ku$; the odd primary case for $\ell$ follows with obvious changes. We start with the $\tau$-Bockstein spectral sequence 
\[E_1=\Z[\tau]\otimes\pi_{*,*}\nu\ku/\tau\implies\pi_{*,*}\nu\ku.\]
    The Omnibus Theorem of \cite[Th.A.1]{burkhahnseng} states that this spectral sequence converges (since $\ku$ is $\F_2$-nilpotent complete as it is connective and implicitly $2$-completed), that there is an isomorphism of $\pi_{*,*}\nu\ku/\tau$ with the $E_2$-page of the ASS of $\ku$, and that differentials in this Bockstein spectral sequence are uniquely determined by differentials in the ASS of $\ku$.\\ 
    
    Part 2 of \Cref{rmk:easyconsequences} then gives us the first isomorphism
    \begin{align*}
        \pi_{*,*}\nu\ku/\tau&\cong\Ext_{\mathcal{A}_*}(\F_2,\A_\ast \square_{\mathcal{E}(1)_\ast} \F_2)\\
        &\cong \Ext_{\mathcal{E}(1)_\ast}(\F_2,\F_2)\\
        &\cong\F_2[h_0,q_1],
    \end{align*}
    where the second line follows from a change-of-rings isomorphism, and the third is the calculation of $\Ext$ over an exterior algebra (see \cite[Th.3.1.16]{greenbook}).\\
    
    The generators $h_0,q_1$ lie in even stems, so the ASS of $\ku$ has no nonzero differentials, and hence the same is true for the $\tau$-Bockstein spectral sequence, and we conclude that the $E_\infty$-page takes the form $\F_2[h_0,q_1,\tau]$.  Choosing a lift $h_0$ defines a surjective map $\Z[h_0,\tau]\to\pi_{0,*}\nu\ku$. It follows that $h_0\tau$ generates the kernel of the $\tau$-reduction map
\[\pi_{0,0}\nu\ku\to\pi_{0,0}\nu\ku/\tau=\F_2\]

Moreover, since the $E_\infty$-page is $\tau$-torsion free, it follows that $\pi_{*,*}\nu\ku$ is $\tau$-torsion free.
    One may form the composite
    \[\Z\to\pi_{0,0}\nu\ku\to\pi_0\ku=\Z\]
    where the first map is defined to send $1$ to $1$, and the second is $\tau$-inversion. The composite is the identity, and hence the second map is a surjection. The second map is then an isomorphism since $\pi_{0,0}\nu\ku$ is $\tau$-torsion free. From this, it follows that the kernel of $\pi_{0,0}\nu\ku\to\pi_{0,0}\nu\ku/\tau$ is precisely $(2)\subset\Z$, and we deduce the relation $h_0\tau=2$, after possibly multiplying $h_0$ by a unit in $\Z$.\\

    Choosing also a lift $q_1$, the resulting map $\Z[h_0,q_1,\tau]/(h_0\tau-2)\to\pi_{*,*}\nu\ku$ is an isomorphism because it induces an isomorphism mod $(\tau)$.    
\end{proof}

This computation leads us to define lifts of $\KU$ and $\L$ to $\Syn_{\F_p}$ with nonvanishing Adams spectral sequences. This is in contrast to the fact that $\nu \KU/\tau$ and $\nu \L/\tau$ are contractible, as $\KU\otimes \F_p$ vanishes.

\begin{mydef}
Set $\KU_{\F_2}:=\nu \ku[q_1^{-1}]$ and for odd primes $\L_{\F_p}:=\nu\ell[q_1^{-1}]$.
\end{mydef}

It is easy to check that $\tau^{-1}\KU_{\F_2}\simeq \KU$ and $\tau^{-1}L_{\F_p}\simeq L$, as $\tau^{-1}$ commutes with colimits and $\tau^{-1}(q_1)=v_1$. Moreover, from \Cref{htpyofcovers}, we see there are equivalences of $\F_p$-synthetic spectra $\tau_{\geq 0}^\uparrow \KU_{\F_p}\simeq \nu\ku$ and similarly that $\tau_{\geq 0}^\uparrow \L_{\F_p}\simeq \nu\ell$.


\subsection{$\nu\ko$ inside $\Syn_{\F_2}$}\label{ftwosynthetickosubsection}

With the bigraded homotopy groups of $\nu\ku$ in our pocket, we can now compute those of $\nu\ko$ courtesy of \Cref{synwood} and the long exact sequence
\begin{equation}\label{syntheticLES}\cdots\to\pi_{n+1,\ast-1}\nu\ku\to\pi_{n-1,\ast-1}\nu\ko\xrightarrow{h_1}\pi_{n,\ast}\nu\ko\to\pi_{n,\ast}\nu\ku\xrightarrow{\partial}\pi_{n-2,\ast}\nu\ko\to\cdots.\end{equation}

By \Cref{filspthm} and \Cref{examplesofsignatures}, this records the data of the (standard) ASS for $\ko$ at the prime $2$. By contrast, the usual recipe to compute the ASS of $\ko$ involves first determining the comodule $H_\ast(\ko;\F_2)$ and then computing $\Ext_{\mathcal{A}(1)_*}(\F_2,\F_2)$. The latter may be done by hand using minimal resolutions or a May spectral sequence, however, we prefer the simpler long exact sequence arguments below.

\begin{figure}[h]\begin{center}
\centering
\makebox[\textwidth]{\includegraphics[trim={3.5cm 17.5cm 3.5cm 4cm},clip,page = 1, scale = 0.8]{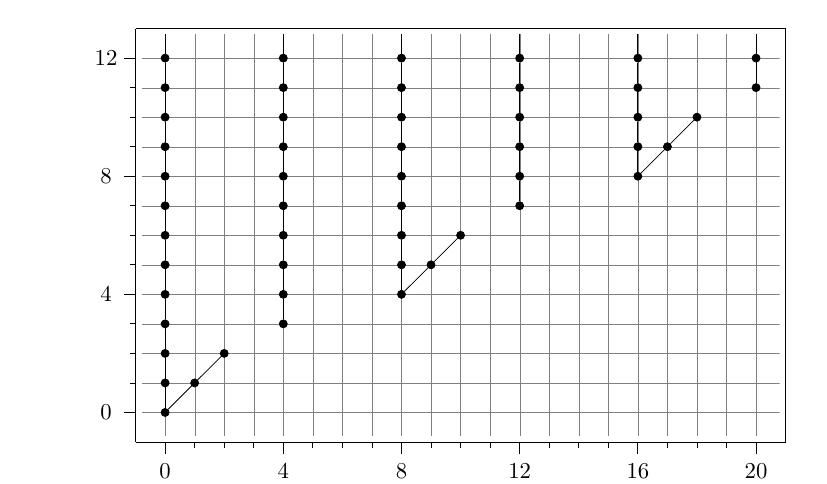}}
\caption{\label{homotopyofkopicture} ASS for $\ko$ at the prime $2$ in the range $0\leq s\leq 20$.}
\end{center}\end{figure}

\begin{figure}[h]\begin{center}
\makebox[\textwidth]{\includegraphics[trim={3.5cm 20.5cm 3.5cm 4cm},clip,page = 1, scale = 0.8]{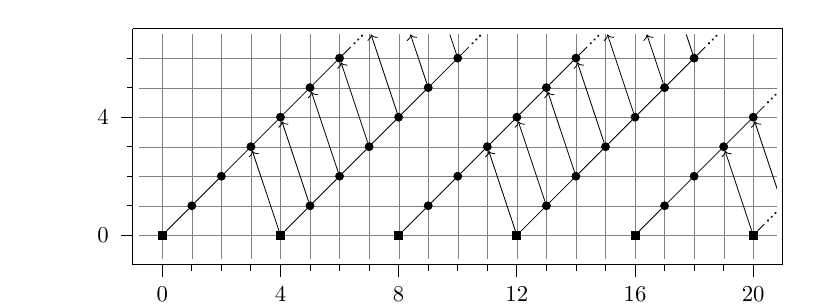}}
\caption{\label{mansskoattwopicture} ANSS for $\ko$ at the prime $2$ in the range $0\leq s\leq 20$.}
\end{center}\end{figure}

\begin{cor}\label{synthetichomotopyofko}
The $\F_2$-synthetic homotopy groups of $\nu\ko$ are given by
\[\pi_{\ast,\ast}\nu\ko=\Z[h_0,h_1,\alpha,\beta,\tau]/(h_0\tau-2,h_0h_1,h_1^3,h_1\alpha,\alpha^2-h_0^2\beta)\]
where $|h_0|=(0,1)$, $|h_1|=(1,1)$, $|\alpha|=(4,3)$, $|\beta|=(8,4)$, and $|\tau|=(0,-1)$. The associated spectral sequence is given as \Cref{homotopyofkopicture}.
\end{cor}

\begin{proof}
For simplicity, write $A$ for the graded ring $\Z[\tau,h_0]/(\tau h_0-2)$. We will show by induction that, when $k\ge0$, there are isomorphisms of $\pi_{0,\ast}\nu\ko$-modules
\[\pi_{8k+r,\ast}\nu\ko\simeq\begin{cases}A\{\be^k\}&r=0\\(A/h_0)\{h_1^r\be^k\}&r=1,2\\A\{\al\be^k\}&r=4\\0&\mathrm{else,}\end{cases}\qquad\qquad |\be|=(8,4),|\al|=(4,3), \]
where the degrees of the generators help us track the $\tau$-degree of these modules, and the class $\beta$ will be defined as a lift of $q_1^4$ along the map $\pi_{8,4}\nu\ko\to\pi_{8,4}\nu\ku$. Simultaneously, we will show that the boundary maps are given (up to a unit) by the formul{\ae}
\[
\partial(q_1^n)=\begin{cases}0&n=4k\\h_0\beta^k&n=4k+1\\h_1^2\beta^k&n=4k+2\\\alpha\beta^k&n=4k+3.\end{cases}
\]
using that, in each case, the claimed value of $\partial(q_1^n)$ is the unique class in the correct bidegree which $\tau$-inverts to the known value of $\partial(u^n)$ in \Cref{cor:woodboundary}. The claimed multiplicative structure then follows from the $h_1$-multiplications and the fact that $\nu\ko\to\nu\ku$ is a ring map.\\

Let us start with the induction step. Suppose that $\pi_{8r,*}\nu\ko\to\pi_{8r,*}\nu\ku$ is an isomorphism for $r\le k$, so that $\pi_{8r,*}\nu\ko\cong A\{\beta^r\}$, where $\beta$ is the unique lift of $q_1^4$. We claim that this determines $\pi_{8k+l,*}\nu\ko$ for $0<l\le 8$, the effect of the boundary map in this range, and that $\pi_{8(k+1),*}\nu\ko\to\pi_{8(k+1),*}\nu\ku$ is also an isomorphism.\\

Indeed, it follows then that $\pi_{8k,*}\nu\ko$ is $\tau$-torsion free, so the boundary map $\pi_{8k+2,*}\nu\ku\to\pi_{8k,*}\nu\ko$ is determined by its effect after inverting $\tau$. It therefore sends $q_1^{4k+1}\mapsto h_0\beta^k$ by \Cref{cor:woodboundary}. Using (\ref{syntheticLES}) along with the fact that $\pi_{\mathrm{odd},\ast}\nu\ku=0$ we see that $\pi_{8k+1,*}\nu\ko=(A/h_0)\{h_1\beta^k\}$, and that $h_1:\pi_{8k+1,\ast}\nu\ko\to \pi_{8k+2,\ast+1}\nu\ko$ is an isomorphism.\\

In particular, we see that $\pi_{8k+2,*}\nu\ko$ is $\tau$-torsion free, so as before we see that $\partial$ sends $q_1^{4k+2}\mapsto h_1^2\beta^k$. By (\ref{syntheticLES}) we see that $\pi_{8k+3,*}\nu\ko=0$, which implies $\pi_{8k+4,*}\nu\ko\to \pi_{8k+4,*}\nu\ku=A$ is an injection with image the ideal generated by $h_0$ (note there is an $A$-module isomorphism $(h_0)\cong \Sigma^{0,1}A$).\\

Then $\pi_{8k+4,*}\nu\ko$ is $\tau$-torsion free, so as before we see that $\partial$ sends $q_1^{4k+3}\mapsto\alpha\beta^k$. Similar arguments to the above now imply that $\pi_{8k+l,*}\nu\ko=0$ for $4<l<8$ and that $\pi_{8k+8,*}\nu\ko\to\pi_{8k+8,*}\nu\ku$ is an isomorphism.\\

This completes the induction step, and for the base case, consider (\ref{syntheticLES}) when $n=0$; the fact that $\ko$ is connective combined with \Cref{tauinversionofcovers} shows that $\nu\ko$ is connective in the vertical $t$-structure, implying that $\pi_{0,*}\nu\ko\to\pi_{0,*}\nu\ku=A$ is an isomorphism by the synthetic Wood cofibre sequence (\ref{syntheticLES}).\\

The computation of the associated spectral sequence (\Cref{homotopyofkopicture}) follows immediately from \Cref{filspthm}; since $\pi_{*,*}\nu\ko$ is $\tau$-torsion free, we have $\pi_{*,*}(\nu\ko/\tau)=(\pi_{*,*}\nu\ko)/\tau$, and there are no nonzero differentials.
\end{proof}

The above calculation also leads to a definition of $\KO_{\F_2}$, an $\F_2$-synthetic $\E_\infty$-ring which produces a nontrivial modified $\F_2$-ASS for $\KO$.

\begin{mydef}
In $\Syn_{\F_2}$, we set $\KO_{\F_2}:=\nu\ko[\beta^{-1}]$.
\end{mydef} 

\begin{cor}\label{taucompletefpstuff}
The $\F_p$-synthetic spectra $\nu \ku,\KU_{\F_p}, \nu\ell, \L_{\F_p},  \nu\ko, \KO_{\F_2}$ are all $\tau$-complete, and applying $\tau^{-1}$ yields the spectra $\ku,\KU,\ell, \L, \ko, \KO$, respectively.
\end{cor}

\begin{proof}
By \cite[Pr.A.13]{burkhahnseng}, the fact that all spectra in sight are implicitly $p$-complete, hence $\F_p$-complete, lead to $\nu\F_p$- and $\tau$-complete synthetic analogues. In the periodic cases, one checks directly that the $\tau$-adic towers are all Mittag-Leffler on bigraded homotopy groups and converge to the bigraded homotopy groups of the synthetic spectrum in question. The claim about $\tau$-inversion follows from the natural equivalence $\tau^{-1}\circ\nu\simeq \id$ and the facts that $\tau^{-1}(q_1)=v_1$ and that $\tau^{-1}$ commutes with colimits.
\end{proof}


\subsection{In $\BP$-synthetic spectra}

These arguments simplify further if we take $E=\BP$. Indeed, in this setting, the ANSSs for $X=\ku$, $\KU$, $\ell$, and $\L$ are incredibly simple, essentially due to the fact that they are complex-oriented and hence have a Thom isomorphism; see part 1 of \Cref{rmk:easyconsequences}. In particular, these ANSSs are concentrated in filtration $0$, hence the bigraded homotopy groups of $\nu X$ in $\Syn_{\BP}$ are $\tau$-torsion free.

\begin{prop}
    There are isomorphisms of bigraded rings
    \[\pi_{\ast,\ast}\nu\ku\simeq \Z[u,\tau],\qquad \pi_{*,*}\nu \KU\simeq \Z[u^{\pm},\tau], \qquad \pi_{*,*}\nu\ell\simeq\Z[v_1,\tau], \qquad \pi_{*,*}\nu L\simeq \Z[v_1^{\pm},\tau]\]
where $|u|=(2,0)$ and $|v_1|=(2p-2,0)$.
\end{prop}

\begin{proof}
    This follows from the exact same arguments as the proof of \Cref{synthetichomotopyofkuandell}, so we will skip some details which are the same as in that proof. In the case of $\nu\ku$, for example, we have a $\tau$-Bockstein spectral sequence of the form
    \[E_1 = \Z[\tau] \otimes \pi_{\ast,\ast}\nu\ku/\tau \implies \pi_{\ast,\ast}\nu \ku.\]
    The discussion above gives us the computation
    \[\pi_{\ast,\ast} \nu\ku/\tau \simeq \Z[u],\qquad |u|=(2,0),\]
    where $u$ detects the class of the same name in $\pi_2 \ku$ in the associated spectral sequence. As there are no differentials in the ANSS for $\ku$, the Omnibus theorem \cite[Th.A.1]{burkhahnseng} tells us there can be no differentials in the above $\tau$-Bockstein spectral sequence. In particular, the $E_1$-page is the $E_\infty$-page. As opposed to the $\F_p$-case, in this world there are no extension problems either, and we arrive at the desired isomorphism. The cases for $\nu\KU$, $\nu\ell$, and $\nu\L$ follow \emph{mutatis mutandis}.
\end{proof}

 As in the previous subsection, we use \Cref{synwood} and the long exact sequence (\ref{syntheticLES}) to compute $\pi_{*,*}\nu\ko$ in $\Syn_\BP$. This determines the entire data of the ANSS of $\ko$.

\begin{prop}\label{homotopykoBP}
The $\BP$-synthetic homotopy groups of $\nu\ko$ are given by
\[\pi_{*,*}\nu\ko=\Z[h_1,\alpha,\beta,\tau]/(2h_1,\tau^2h_1^3,h_1\alpha,\alpha^2-4\beta)\]
where $|h_1|=(1,1)$, $|\alpha|=(4,0)$, $|\beta|=(8,0)$, and $|\tau|=(0,-1)$,
The associated spectral sequence is given as \Cref{mansskoattwopicture}.
\end{prop}

Notice the similarities to \cite[Th.4.10]{isasksenshkembi}, the only different being that their $\tau$ acts as our $\tau^2$---this is the difference between $\BP$-synthetic spectra and even $\BP$-synthetic spectra; see \cite[5.2]{syntheticspectra}.

\begin{proof}
The argument closely resembles the $\Syn_{\F_2}$-case, so we will be brief. The crucial point is to compute the effect of the boundary map in (\ref{syntheticLES}) in $\Syn_\BP$; the computation then follows as usual by computing the kernel and cokernel. Note that the boundary maps $\pi_{2k,*}\nu\ku\to\pi_{2k-2,*}\nu\ko$ are determined by their effect on $u^k\in\pi_{2k,0}\nu\ku$ as maps of $\Z[\tau]$-modules. The boundary maps are therefore determined by their effect after inverting $\tau$, as $\pi_{n,0}\nu\ko$ is always $\tau$-torsion free. Indeed, a nonzero $\tau$-torsion class in $\pi_{n,0}\nu\ko$ would imply the existence of a nonzero differential in $\mathrm{ANSS}(\ko)$ with target in bidegree $(n,0)$; this is impossible as standard ANSSs are concentrated in nonnegative filtrations.
In particular, we find that (up to a unit)
\[
\partial(u^n)=\begin{cases}0&n=4k\\2\beta^k&n=4k+1\\\tau^2h_1^2\beta^k&n=4k+2\\\alpha\beta^k&n=4k+3\end{cases}
\]
as, in each case, the claimed value of $\partial(u^n)$ is the unique class in the correct bidegree which $\tau$-inverts to the known value of $\partial(u^n)$ in \Cref{cor:connectivewoodstheorem}.
\end{proof}

Either by inverting $\beta$ or using a similar argument as above with $\nu\KO$, one may also compute $\pi_{*,*}\nu\KO$.\\

In the category $\Syn_\BP$, it is noteworthy that one may construct synthetic lifts of $\KO$ and $\ko$ directly from $\nu \KU$, and these turn out to recover the synthetic analogues $\nu\KO$ and $\nu\ko$. The rest of this section is not necessary for the main purposes of this article, we only include it to provide further intuition and warnings for both the reader and ourselves. Recall that $\KO_{\BP} = (\nu \KU)^{hC_2}$ from \Cref{definitionofko}.

\begin{mydef}\label{defkobp}
Let $L$ be the line $y=x$. In $\Syn_{\BP}$, we set $\ko_{\BP} = \tau_{\ge0}^L \KO_{\BP}$.
\end{mydef}

\begin{prop}\label{agreementinbpcase}
There are equivalences of $\BP$-synthetic spectra $\Syn_\BP$, $\nu\KO\simeq\KO_\BP$ and $\nu\ko\simeq\ko_\BP$.
\end{prop}

\begin{proof}
There are natural assembly maps $\nu\KO\to\KO_\BP$ and $\nu\ko\to\ko_\BP$, and it suffices to show they induce isomorphisms on bigraded homotopy groups. We may compute the bigraded homotopy groups of $\KO_{\BP}$ (and thereby of $\ko_{\BP}$) by running the homotopy fixed point spectral sequence (HFPSS) internal to $\Syn_{\BP}$. This HFPSS has signature
\[E_2=H^s(C_2;\pi_{n,f}\nu\KU)\implies\pi_{n-s,f+s}\KO_\BP\]
drawn in Adams trigrading $(n-s,s+f,f)$, the differential $d_r$ has signature $(-1,+1,+r)$.\\

The synthetic homotopy groups $\pi_{*,*}\nu\KU=\Z[u,\tau]$ are $\tau$-torsion free, so it follows that the generator of $C_2$ fixes $\tau$ and sends $u\mapsto-u$, as this is true after inverting $\tau$. Standard group cohomology calculations then give
 \[E_2=\Z[h_1,u^{\pm 2},\tau]/(2h_1)\]
This is, in particular, $\tau$-torsion free, and one recovers the classical HFPSS for $\KU$ upon $\tau$-inversion, since applying $\tau^{-1}$ to the cosimplicial object
\[\nu\KU\implies F({C_2}_+,\nu\KU)\Rrightarrow F({C^{\times2}_2}_+,\nu\KU)\cdots\]
gives the cobar resolution for $\KU$. In the classical HFPSS for $\KU$, there is a differential $d_3(u^2)=\eta^3$, which may be deduced from the fact that $\eta^4=0\in\pi_4\Sph$, which follows from our basic assumptions about $\pi_\ast \Sph$ in low degrees; see \Cref{sec:basicassumptions}. Therefore $\tau^{-1}(d_3(u^2))=\eta^3$, and synthetically, $\tau^2h_1^3$ is the unique class in the tridegree of $d_3(u^2)$ that $\tau$-inverts to $\eta^3$, so we deduce $d_3(u^2)=\tau^2h_1^3$. The HFPSS collapses on the $E_4$ page for degree reasons, and the $E_\infty$ page is isomorphic to the bigraded homotopy groups of $\nu\KO$ as in \Cref{homotopykoBP}. There are no nontrivial extensions as each relation holds in or above the highest nonzero filtration in its respective stem on the $E_\infty$ page of the HFPSS.
\end{proof}


\begin{warn}\label{koF2nogo}
One may define a $\tau$-complete $\F_2$ synthetic lift of $\ko$ in a similar way to \Cref{defkobp} by setting
\[\ko_{\F_2}:=\tau_{\ge0}^\uparrow(\nu\ku)^{hC_2}\]
This is, however, not the same as $\nu\ko$. Indeed, from the HFPSS internal to $\Syn_{\F_2}$ with signature
\[E_2=H^s(C_2;\pi_{n,f}\ku_{\F_2})\implies \pi_{n-s,f+s}(\KU_{\F_2})^{hC_2}.\]
the nonzero class in $H^1(C_2;\pi_{2,1}\ku_{\F_2})$ detects a class $x\in\pi_{1,2}(\nu\ku)^{hC_2}$ with $h_0x\neq0$ and $2x=0$. This implies that $h_0x$ is a nonzero $\tau$-torsion class. However, it follows from \Cref{synthetichomotopyofko} that $\pi_{\ast,\ast}\nu\ko$ is $\tau$-torsion free, hence $\ko_{\F_2}$ cannot be equivalent to $\nu\ko$. In fact, the MASS associated with $\ko_{\F_2}$ turns out to be very similar to the spectral sequence of Hill--Lawson converging to the sphere spectrum \cite{hill2021ekpushouts}. We intend to return to this and similar examples in future work.
\end{warn}

The following is simply \cite[Pr.A.13]{burkhahnseng} and the natural equivalence $\tau^{-1}\circ\nu\simeq \id$.

\begin{cor}\label{taucompetenessofbpstuff}
The $\BP$-synthetic spectra $\nu \ku,\nu\KU,\nu\ell,\nu \L, \nu\ko, \nu \KO$ are all $\tau$-complete, and applying $\tau^{-1}$ gives the spectra $\ku,\KU,\KO,\ko,\ell,\L$, respectively.
\end{cor}


\section{Two synthetic image-of-$J$ spectra}\label{imageofjconstructionsection}

Classically, one defines the connective image-of-$J$ spectrum $\j$ at the prime $2$ as follows: first, one shows that the map of spectra $\psi^3-1\colon \ko\to \ko$ factors through $\tau_{\geq 4}\ko$, and then one defines $\j$ as the fibre of this factorisation. To construct our synthetic lifts of $\j$, we need to show that the map of synthetic spectra $\psi^3-1\colon \nu\ko\to \nu\ko$ factors through a connective cover with respect to our vertical $t$-structure.


\subsection{Construction of synthetic lifts of $\j$}

To construct these lifts, we need stable Adams operations. Recall that everything is implicitly $p$-completed.

\begin{mydef}
    Let $E$ be either $\F_p$ or $\BP$. If $p=2$, we write $\psi^3\colon \nu\ko \to \nu\ko$ in $\Syn_E$ for $\nu \psi^3$. If $p$ is odd, we write $\psi^{p+1}\colon \nu \ell \to \nu\ell$ in $\Syn_E$ for $\nu \psi^{p+1}$.
\end{mydef}

\begin{prop}\label{factorisationthroughcovers}
    Let $E$ be either $\F_p$ or $\BP$. Then there is a factorization of the maps
    \[\psi^3-1:\nu\ko\to \nu\ko\qquad \qquad \psi^{p+1}-1\colon \nu\ell\to \nu\ell\]
    through $\tau_{\ge 4}^\uparrow\nu\ko$ and $\tau_{\geq 2p-2}^\uparrow\nu\ell$ at $p=2$ and odd primes, respectively.
\end{prop}

\begin{proof}
There is a fibre sequence of synthetic spectra $\tau_{\geq 3}^\uparrow\nu\ko\to \nu\ko\to\tau_{\leq 2}^\uparrow\nu\ko$. Notice that the composite
    \[\nu\ko\xrightarrow{\psi^3-1}\nu\ko\to\tau_{\le 2}^\uparrow\nu\ko,\]
    or its adjunct $\tau_{\le 2}^\uparrow\nu\ko\xrightarrow{\tau_{\le 2}^\uparrow(\psi^3-1)} \tau_{\le 2}^\uparrow\nu\ko$, is null. Indeed, as the truncation of the unit $\tau_{\le 2}^\uparrow\1\to\tau_{\le 2}^\uparrow\nu\ko$ is an equivalence for both $E=\BP$ and $\F_p$, then this adjunct is equivalent to a map $\tau_{\le 2}^\uparrow\1\to \tau_{\le 2}^\uparrow\nu\ko$, which is in turn adjoint to the map
    \[(\psi^3-1)(1) \colon \1\to \tau_{\le 2}^\uparrow\nu\ko.\]
    This map is null as it represents zero in $\pi_{0,0}\tau_{\leq 2}^\uparrow \nu\ko$. This shows the desired map factors through $\tau_{\geq 3}^\uparrow \nu\ko$. As $\pi_{3,\ast}\nu\ko=0$ for $E=\F_2$, see \Cref{homotopyofkopicture}, it follows that $\tau_{\geq 3}^\uparrow \nu\ko=\tau_{\geq 4}^\uparrow \nu\ko$, so the $E=\F_2$-case is done. Similarly, for an odd prime $p$, we have $\pi_{2p-3,\ast}\nu\ell=0$ for both $E=\F_p$ and $\BP$, so these cases are also covered. For $E=\BP$ at $p=2$, we have the following solid diagram of synthetic spectra
    \[\begin{tikzcd}
    	&	{\tau_{\geq 4}^\uparrow \nu\ko}\ar[d]	&	\\
	{\nu\ko}\ar[r, "\psi^3-1"]\ar[ru, dashed]	&	{\tau_{\geq 3}^\uparrow \nu\ko}\ar[r]	&	{\Sigma^{3,3}\F_2[\tau]/\tau^2}
    \end{tikzcd}\]
where the composite of the vertical and the right map is a fibre sequence; see \Cref{mansskoattwopicture}. To obtain the dashed arrow, we want to show that the lower composite is null. This follows from another series of formal manipulations. This composite precomposed with the unit $\1\to \nu\ko$ vanishes for degree reasons, so the above composite factors through the cofibre of the unit, denoted by $\nu\ko/\1$. As this unit induces an isomorphism on $\pi_{n,\ast}$ for $n\leq 2$ and a surjection on $\pi_{3,\ast}$, we see that $\tau_{\leq 3}^\uparrow (\nu\ko/\1)=0$. In particular, we see the factorisation of the above composite is adjunct to the zero map
\[0=\tau_{\leq 3}^\uparrow \nu\ko/\1\xrightarrow{\tau_{\leq 3}^\uparrow (\psi^3-1)} \Sigma^{3,3}\F_2[\tau]/\tau^2.\qedhere\]
\end{proof}

Our synthetic lifts of $\j$ are then defined just as in the classical case, now using the vertical $t$-structure of \Cref{tstructuresexist}.

\begin{mydef}\label{imageofjsyntheticdefinition}
    Let $E$ be either $\F_p$ or $\BP$. By \Cref{factorisationthroughcovers}, we have maps of synthetic spectra
    \[\psi^3-1\colon \nu\ko\to \tau_{\geq 4}^\uparrow \nu\ko\qquad \qquad \psi^{p+1}-1\colon \nu\ell\to \tau_{\geq 2p-2}^\uparrow\nu\ell\]
    where $p$ is an odd prime. Define $\j_E$ as the fibre of the above maps. As the above maps are zero on $\pi_{0,0}$, we obtain a map $\1\to \j_E$ which we call the \emph{unit map}; in \Cref{einftystructuresonj} we will see these synthetic spectra have preferred $\E_\infty$-structures.
\end{mydef}

First, let us show $\j_E$ are indeed lifts of the spectra $\j$.

\begin{prop}\label{jearelifts}
    For all primes $p$, the synthetic spectra $\j_{\F_p}$ and $\j_{\BP}$ are $\tau$-complete synthetic lifts of the connective image-of-$J$ spectrum $\j$.
\end{prop}

\begin{proof}
    By \Cref{taucompetenessofbpstuff,taucompletefpstuff,coverstaucomplete}, and the fact that $\tau$-inversion is exact, we see that $\j_E$ is $\tau$-complete, and the fibre sequences of synthetic spectra defining $\j_E$ are mapped precisely to the fibre sequences of spectra defining $\j$. 
\end{proof}


\subsection{Synthetic homotopy groups of $\j_E$}

To highlight the usefulness of our definition of $\j_E$, let us calculate $\pi_{\ast,\ast}\j_E$---by \Cref{jearelifts}, these synthetic homotopy groups are the signature of a modified $E$-ASS for the spectrum $\j$. The proofs of all of the following calculations are the same, so we will only prove the first two statements at the prime $2$.

\begin{prop}\label{evenassjcalculation}
    Write $A$ for the graded ring $\Z[\tau,h_0]/(\tau h_0-2)$. For $k\in\Z$ and $0\leq r\leq 7$, there are isomorphisms of graded $A$-modules
    \[
    \pi_{8k+r,\ast}(\j_{\F_2})\simeq \left\{\begin{array}{lll}
	A\{1\}							&	k=r=0		\\
	(A/h_0)\{P^{k-1}(h_1h_3)\}				&	k>0, r=0		\\
	(A/h_0)\{h_1^r\}						&	k=0, r=1,2		\\
	(A/h_0)\{P^{k-1}(h_1^2h_3), P^k(h_1)\}	&	k>0, r=1		\\
	(A/h_0)\{P^k(h_1^2)\}					&	k>0, r=2		\\
	(A/(\tau h_0)^3)\{P^k(h_2)\}					&	k\geq 0, r=3 	\\
	(A/(\tau h_0)^{\ord_2(k+1)+4})\{\al_k\}			&	k\geq 0, r=7	\\
	0								&	\text{otherwise}	&	
    \end{array}\right.\]
	with bidegrees on multiplicative generators given by $|h_i|=(2^i-1,1)$, $|\alpha_k|=(7+8k,5+4k)$, and $P^k(x)=|x|+(8k,4k)$.
\end{prop}

The notation $P^k(-)$ is supposed to suggest the existence of a ``synthetic Adams periodicity operator''; see \Cref{periodicityoperatorfailsnot}.

\begin{proof}
	The computation follows by determining the effect of the map $\psi^3-1:\nu\ko\to\tau_{\ge4}^\uparrow\nu\ko$ on bigraded homotopy groups. This is easy, as the source and target are $\tau$-torsion free, and we know the effect after inverting $\tau$; we see that the map is zero in $\pi_{8k+r,\ast}$ when $r\neq0,4$, and it is given by multiplication by $8=\tau^3h_0^3$ when $r=4$ and multiplication by $2^{\ord_2(k)+4}=(\tau h_0)^{\ord_2(k)+4}$ when $r=0$ and $k>0$, each of these holding up to multiplication by an odd number, in particular a $2$-adic unit. This finishes our calculations up to extensions---the rest of the proof is simply a resolution of a classical extension problem calculating $\pi_d\j$ for $d\equiv_8 1$. \\

For degree reasons, the only possible nontrivial extensions in the long exact sequence computing $\pi_{\ast,\ast}\j_{\F_2}$ appear in stems of the form $8k+1$. For a contradiction, suppose this extension was not trivial so that after inverting $\tau$ we have $\pi_{8k+1}\j\simeq \Z/4\Z$. A quick calculation from the long exact sequence on homotopy groups associated with $\Sph/2$ would show that $\pi_{8k+1}\j/2$ has $4$-elements. However, the long exact sequence on homotopy groups associated with the cofibre sequence
\[\j/2\to \ko/2\xrightarrow{(\psi^3-1)/2} \tau_{\geq 4}\ko/2\]
produces the exact sequence
\begin{equation}\label{modtwoles}\pi_{8k+2}\ko/2\xrightarrow{\psi^3-1} \pi_{8k+2}\tau_{\geq 4}\ko/2\to \pi_{8k+1}\j/2\to \pi_{8k+1}\ko\xrightarrow{\psi^3-1=0}\pi_{8k+1}\tau_{\geq 4}\ko.\end{equation}
If there is an isomorphism $\pi_{8k+2}\ko/2\simeq\Z/4\Z$ and $\psi^3-1$ acts trivially on this group, then we are done. Indeed, as $\pi_{8k+1}\ko/2\simeq \Z/2\Z$, which easily follows from the long exact sequence on homotopy groups, we see that $\pi_{8k+1}\j/2$ must have $8$ elements from (\ref{modtwoles}), hence $\pi_{8k+1}\j$ must be isomorphic to $\Z/2\Z\oplus\Z/2\Z$, as desired. To see that $\pi_{8k+2}\ko/2\simeq \Z/4\Z$, we use $8$-fold Bott periodicity, combined with the facts that $\Sph\to \ko$ induces an isomorphism on $\pi_{s}$ for $s=1,2$ and the classical calculation $\pi_2\Sph/2\simeq \Z/4\Z$; see \Cref{anssmooreattwopicture}. To see that $\psi^3-1$ acts trivially on $\pi_{8k+2}\ko/2$, it suffices to show that $\psi^3$ acts by the identity. To see this, take a generator $\bar{\be}^{k}x$ of this group, where $x$ is a generator of $\pi_2\ko/2\simeq\Z/4\Z$. Then $\psi^3(\bar{\be}^k x)=\psi^3(\bar{\be})^k\psi^3(x)$ by multiplicativity, and as $\psi^3(\bar{\be})=\bar{\be}$, we are reduced to the $\pi_2$-case again. However, $\psi^3$ acts trivially on $\pi_2\ko/2$ as the unit $\Sph\to \ko$ induces an isomorphism on $\pi_s$ for $s=1,2$, and $\psi^k$ commutes with this unit as a map of spectra, hence its reduction mod $2$ also commutes with this isomorphism $\pi_2\Sph/2\simeq \pi_2\ko/2$.
\end{proof}

In the classical ASS for $\j$, there are no such extension problems in $\pi_{8k+1}\j$ to consider for degree reasons.

\begin{figure}[h]\begin{center}
\makebox[\textwidth]{\includegraphics[trim={2.5cm 15cm 2.5cm 4cm},clip,page = 2, scale = 0.8]{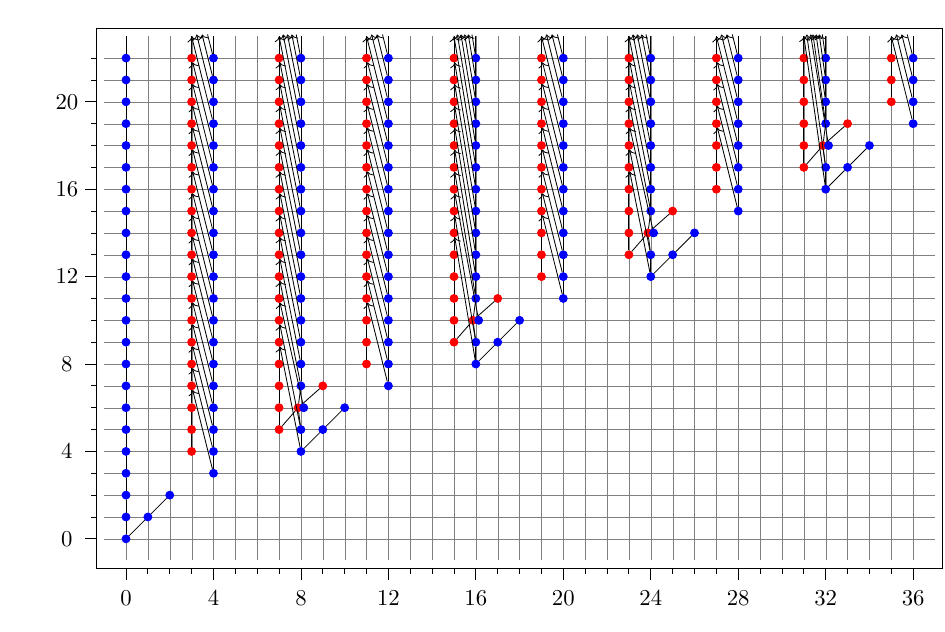}}
\caption{\label{massjattwopictureetwo} MASS for $\j$ at the prime $2$ (from $\j_{\F_2}$) in the range $0\leq s\leq 36$ with all differentials.}
\end{center}\end{figure}

\begin{cor} 
The MASS for $\j$ associated with the $\F_2$-synthetic spectrum $\j_{\F_2}$ has $E_2$-page given by $\pi_{\ast,\ast}(\nu\ko/\tau)\oplus\pi_{\ast,\ast}(\tau_{\ge4}^\uparrow(\nu\ko)/\tau)[-1]$. There are $d_4$'s leaving stems of the form $8k+4$ and $d_{\ord_2(k)+5}$'s leaving the $8k$ stem, as shown in \Cref{massjattwopictureetwo}. Moreover, there are exotic $\eta$-extensions indicated by the orange lines of slope $4$ in \Cref{massjattwopictureetwo}.
\end{cor}

\begin{proof}
We see from \Cref{evenassjcalculation} that $\psi^3-1$ induces the zero map on $\pi_{\ast,\ast}(-/\tau)$ from which the claim about the $E_2$ page follows. The differentials follow immediately from the $\tau$-torsion described in \Cref{evenassjcalculation}. For the $\eta$-extensions, we need to use the modified ANSS for $\j$, which is easy to describe from \Cref{evenanssjcalculation}; also see \Cref{manssjattwopictureetwo}. In particular, in this modified ANSS for $\j$ is clear that in $\pi_\ast\j$ we have $\eta P^k(\eta^2)=4 P^k(\nu)$, where $P^k(\eta^2)$ is the class in $\pi_{8k+2}\j$ detected by $P^k(h_1^2)$ and $P^k(\nu)$ the class in $\pi_{8k+3}\j$ detected by $P^k(h_2)$. This necessitates the exotic $\eta$-extensions of \Cref{massjattwopictureetwo}.
\end{proof}

Then the $E=\BP$-case at $p=2$. Because of the $\eta$-towers in $\pi_{*,*}\nu\ko$, it is difficult to describe $\pi_{*,*}\j_{\BP}$ stem-by-stem. We make the following statement instead:

\begin{prop}\label{evenanssjcalculation}
    Let $R$ be the ring $\Z_2[\tau,\eta]/(2\eta,\tau^2\eta^4)$, and for $k\ge0$, let $B_k$ be the ring $\Z/2^{\ord_2(k+1)+4}[\tau,\eta]/(2\eta,\tau^2\eta^3)$, where $|\eta|=(1,1)$ and $|\tau|=(0,-1)$. Then there is an isomorphism of $\Z_2[\tau]$-modules
    \[\pi_{*,*}\j_\BP=R\{1,\theta_{3,5},\nu\}/(2\theta_{3,5},\tau^2\eta^3\theta_{3,5},\eta\nu,4\nu-\tau^2\eta^3)\oplus\bigoplus\limits_{k\ge0}B_k\{\alpha_k,x_k,y_k\}/(2x_k,\eta y_k,4y_k-\tau^2\eta^2x_k)\]
    where $|\theta_{3,5}|=(3,5)$, $|\nu|=(3,1)$, $|\alpha_k|=(8k+7,1)$, $|x_k|=(8k+9,1)$, $|y_k|=(8k+11,1)$.
\end{prop}

\begin{proof}
We again proceed by calculating the long exact sequence in bigraded homotopy groups associated with the fibre sequence defining $\j_{\BP}$. Here $\pi_{*,*}\ko_{\BP}$ is no longer $\tau$-torsion free, so it is easier to calculate the effect of the ring map $\psi^3:\ko_{\BP}\to\ko_{\BP}$ by checking on the generators, use that $\tau_{\ge4}^\uparrow\ko_{\BP}\to\ko_{\BP}$ induces an injection on bigraded homotopy groups, and then calculate $\psi^3-1$. Since $\eta$ and $\tau$ are in the image of the unit map to $\ko_{\BP}$, we see that $\psi^3(\eta)=\eta$ and $\psi^3(\tau)=\tau$. It follows from \Cref{homotopykoBP} that $\pi_{4,0}\ko_{\BP}$ and $\pi_{8,0}\ko_{\BP}$ are $\tau$-torsion free, so we see that $\psi^3(\alpha)=9\alpha$ and $\psi^3(\beta)=81\beta$ since this holds after inverting $\tau$.\\

The calculation now follows from \Cref{homotopykoBP}; it is best to think of the calculation of these homotopy groups in terms of the spectral sequence they correspond to, which is depicted in \Cref{manssjattwopictureetwo} (here the red dots come from $\pi_{*,*}(\tau_{\ge 4}^\uparrow(\nu\ko)/\tau)$ via the boundary map, and the blue dots are sent to the corresponding classes in $\pi_{*,*}(\nu\ko/\tau)$). The only nontrivial point is the relations $4\nu=\tau^2\eta^3$ and $4y_k=\tau^2\eta^2x_k$, which follow from the fact that the corresponding stem must $\tau$-invert to known homotopy groups of $\j$ which are isomorphic to $\Z/8$.
\end{proof}

As per usual, the odd prime case for $E=\F_p$ and $\BP$ are much easier; we omit proofs.

\begin{prop}\label{oddjcalculation}
    Let $p$ be an odd prime and write $q=2p-2$. Write $C$ for the graded ring $\Z[\tau,h_0]/(h_0\tau-p)$ and $D$ for the graded ring $\Z[\tau]$. For $n\in \Z$ we have isomorphisms of $C$- and $D$-modules
  \[
    \pi_{n,\ast}(\j_{\F_p})\simeq \left\{\begin{array}{lll}
	C\{1\}							&	n=0		\\
	C/h_0^{\ord_p(n)+2}\{\al_{n/\F_p}\}		&	n\geq 0, n\equiv -1 \text{ modulo }q	\\
	0								&	\text{otherwise}	&	
    \end{array}\right.\]
      \[
    \pi_{n,\ast}(\j_{\BP})\simeq \left\{\begin{array}{lll}
	D\{1\}							&	n=0		\\
	D/p^{\ord_p(n)+2}\{\al_{n/\BP}\}			&	n\geq 0, n\equiv -1 \text{ modulo }q	\\
	0								&	\text{otherwise}	&	
    \end{array}\right.\]
    respectively, where $|\al_{n/\F_p}|=(nq-1,\frac{nq}{2}+1)$ and $|\alpha_{n/\BP}|=(nq-1,1)$.
\end{prop}

Note that the calculation of the synthetic homotopy groups of $\j_{\F_p}$ is only slightly more complicated for $p=2$, in stark contrast to the classical situation, where the ASS for $\j$ at odd primes is easy to produce, but that for $p=2$ has only recently appeared in the literature; see \cite{imageofjbrunerrognes}. Comparing our calculation above to these calculations of Bruner--Rognes, we also see that $\j_{\F_p}$ is \textbf{not} equivalent to $\nu \j$.


\subsection{Multiplicative structure on $\j_E$}

Using these synthetic homotopy groups, we can show $\j_E$ can be given an $\E_\infty$-structure. The linear $t$-structures used in this section are compatible with the symmetric monoidal structure on $\Syn_E = \Syn$ by \Cref{tstructuresexist}. In particular, both the the connective cover functor $\tau_{\geq 0}^L\colon \Syn \to \Syn_{\geq 0}^L$ and the restricted truncation functor $\tau_{\leq 0}^L\colon \Syn_{\geq 0}^L \to \Syn^L_{\leq n}\cap \Syn^L_{\geq 0} = \Syn^L_{[0,n]}$ are both lax monoidal, hence send $\E_\infty$-objects to $\E_\infty$-objects. To see this, suppose $L$ and $E$ are as in \Cref{tstructuresexist}, such that the associated linear $t$-structure on $\Syn_E=\Syn$ is monoidal. Then $\Syn^L_{\geq0}$ is closed under finite tensoring and contains the unit, hence it inherits a symmetric monoidal structure from $\Syn$ such that the inclusion $\Syn_{\geq0}^L\to \Syn$ is strong symmetric monoidal. The right adjoint, $\tau_{\geq 0}^L\colon \Syn \to \Syn_{\geq 0}^L$ is then lax symmetric monoidal by \cite[Pr.2.5.5.1]{sagname}. On the other hand, \cite[Ex.2.2.1.10]{haname} shows us that the restricted truncation functor above is also lax symmetric monoidal.

\begin{prop}\label{einftystructuresonj}
    Let $E$ be either $\F_p$ or $\BP$. Write $\j'_E$ for the synthetic $\E_\infty$-ring defined by taking the homotopy fixed points of the $\Z$-action on $\nu\ku$ generated by $\psi^3$ at the prime $2$, and the $\Z$-action on $\nu\ell$ generated by $\psi^{p+1}$ at odd primes. Then there is a Cartesian diagram of synthetic spectra
  \[\begin{tikzcd}
        {\j_E}\ar[r]\ar[d]  &   {\tau_{\geq 0}^\uparrow \j'_E}\ar[d] \\
        {\tau_{\leq 2}^\uparrow\1}\ar[r]  &   {\tau_{\leq 2}^\uparrow \tau_{\geq 0}^\uparrow \j'_E.}
    \end{tikzcd}\]    
    In particular, the lower-horizontal and right-vertical maps are maps of synthetic $\E_\infty$-rings, so $\j_E$ can be endowed with the structure of a synthetic $\E_\infty$-ring.
\end{prop}

\begin{proof}
Suppose that $p=2$ for simplicity. One calculates the synthetic homotopy groups of $\j_E'$ as it can be written as the fibre of the map $\psi^3-1\colon \nu\ko\to \nu\ko$. From this one easily computes the synthetic homotopy groups of the vertical connective cover. One can then calculate the synthetic homotopy groups of the actual pullback, which we will denote by $P$, using a Mayer--Vietoris sequence. In total, we find that the synthetic homotopy groups of $P$ match those of $\j_E$. Moreover, at the prime $2$, as the composite
    \[P\to \j'_E\to \nu\ko\xrightarrow{\psi^3-1}\tau_{\geq 4}^{\uparrow}\nu\ko\]
    is null by construction, we obtain a map of synthetic spectra $\j_E\to P$ that recognises the above isomorphism on synthetic homotopy groups. The argument at odd primes follows \emph{mutatis mutandis}.
\end{proof}


\subsection{Periodic variants}\label{periodicvariantsofJsubsection}

Just as in \Cref{ktheoryconstructionsection}, we can also define periodic versions of $\j_E$, denoted by $\J_E$, which behave both as synthetic lifts of the classical periodic image-of-$J$ spectrum and produce a nonzero modified $E$-ASS.

\begin{mydef}
Let $p$ be a prime and $E=\F_p$ or $\BP$. Then we define $\J_E$ as the fibre of
\[\KO_E\xrightarrow{\psi^3-1} \KO_E \qquad\qquad \L_E\xrightarrow{\psi^{p+1}-1}\L_E\]
for $p=2$ and odd primes, respectively.
\end{mydef}

The following two statements are proven just as \Cref{jearelifts,evenassjcalculation}; we omit the proofs. Calculations of $\pi_{\ast,\ast}\J_{\F_p}$ at other primes follow similarly.

\begin{prop}
The synthetic spectra $\J_E$ are $\tau$-complete synthetic lifts of $J$.
\end{prop}

\begin{prop}
At the prime $2$, the synthetic homotopy groups of $\J_{\F_2}$ and $\J_\BP$ are given by \Cref{assperiodicjatinftytwoetwopicture,assperiodicjattwoeinftypicture} and \Cref{anssperiodicjattwoetwopicture,anssperiodicjeinftypicture}, respectively.
\end{prop}

Of course, the odd primary cases are simpler, so we will omit them here.\\

As $\J$ is equivalent to the $\K(1)$-local sphere, we see that \Cref{assperiodicjatinftytwoetwopicture,assperiodicjattwoeinftypicture} provide a modified $\F_p$-ASS for $\Sph_{\K(1)}$. A similar conclusion can be made about a modified ANSS for $\Sph_{\K(1)}$ using the synthetic homotopy groups of $\J_{\BP}$. In this case, however, we claim that $\nu\J\simeq \J_{\BP}$, so this simply recovers the usual ANSS for the $\K(1)$-local sphere. This calculation is well-known, and a modern reference and diagram can be found in \cite[Fig.11]{itamarknbrauer}.


\section{Adams $v_1$ self-maps and detection for $\ko$}\label{selfmapssection}

This section follows from the classical work of Adams \cite{adamsjx}, however, the our proofs below will both act as a warm-up for \Cref{detectionresultsforj} and allow us to solidify our notation.


\subsection{Construction of $v_1$ self-maps}

To show that the unit map $\Sph\to \j$ is surjective on $\pi_{8k-1}$, we will use Adams' $v_1$ self-map $v_1^4\colon \Sph/2[8]\to \Sph/2$; see \cite{adamsjx}. Let us reprove this classical result with the tools we have at hand---this will help us set up the following detection statements as well. Let us write $u$ for the generator of $\pi_2 \ku/2$.

\begin{theorem}[{\cite{adamsjx}}]\label{adamsperiodicprimetwo}
    There is a map of spectra $v_1^4\colon \Sph/2[8]\to \Sph/2$ inducing multiplication by $u^4$ on $\ku$-homology.
\end{theorem}

Recall from \Cref{sec:basicassumptions} that we take for granted the ASS and ANSS for $\Sph$ in stems $\leq 8$.

\begin{proof}
As (the rotation of) the cofibre sequence defining $\Sph/2$ induces a short exact sequence on $\F_2$-homology, then according to \cite[Lm.4.23]{syntheticspectra}, one has a cofibre sequence in $\Syn_{\F_2}$
\[\Sph^{0,1}\xrightarrow{h_0}\nu(\Sph)\to\nu(\Sph/2)\to\nu(\Sph^1)\xrightarrow{h_0}\Sph^{1,0}.\]
From this, one can easily compute the Adams chart for $\Sph/2$ in stems $\le 8$, as a mod $\tau$-computation reveals there cannot be any differentials in this range; see \Cref{assmooreattwopicture}. Similarly,  one has a cofibre sequence in $\Syn_\BP$
\[\1\xrightarrow{2}\1\to\nu(\Sph/2)\to\Sigma^{1,-1}\1\xrightarrow{2}\Sigma^{1,-1}\1\]
which yields an AN-$E_2$-page for $\Sph/2$ in stems $\le 8$. There is a collection of $d_3$-differentials stemming from the key $d_3$-differential in the ANSS for $\Sph$ hitting $h_1^4$. As a module over the ANSS for $\Sph$, the ANSS for $\Sph/2$ is determined by this differential; see \Cref{anssmooreattwopicture,anssmooreattwopictureeinfty}. In both the ASS and ANSS, the exotic $2$-extension on their abutments comes from our basic assumption that $\pi_2 \Sph/2\simeq \Z/4\Z$.\\

We claim that $2\cdot \pi_8(\Sph/2)=0$. Indeed, looking at the long exact sequence
\[\pi_8 \Sph \xrightarrow{2} \pi_8 \Sph \to \pi_8 \Sph/2 \xrightarrow{\partial} \pi_7 \Sph \xrightarrow{2} \pi_7 \Sph,\]
we see the cokernel of the first map is $\F_2\{\eta\sigma, \varepsilon\}$, and the kernel of the last map is generated by $8\sigma$. From the ASS for $\Sph/2$, so \Cref{assmooreattwopicture}, we see that the lift of $h_0^3h_3$ in $\pi_{8,4}\1/(2,\tau)$ which represents a choice of lift of $8\sigma$ has filtration $4$, a higher filtration than that of the image of $h_1h_3$ and $c$ in $\pi_{8,2}\1/(2,\tau)$ and $\pi_{8,3}\1/(2,\tau)$, respectively. In particular, this choice of lift of $h_0^3h_3$ is $2$-torsion, and $\pi_8 \Sph/2 \simeq (\F_2)^3$ as there are no differentials for degree reasons. It follows that any map $\Sph^8\to \Sph/2$ factors over $\Sph/2[8]$.\\

It remains to use the ANSS for $\Sph/2$ (\Cref{anssmooreattwopicture,anssmooreattwopictureeinfty}) to show that there is some class $v_1^4\in\pi_8\Sph/2$ which is sent to $u^4\in\pi_8\ku/2$. There is such a class on the $E_2$-page of the ANSS for $\Sph/2$ from the calculation
\[\pi_{\ast,0}\1/(\tau,2) \simeq \Ext_{\BP_*\BP}^{0,*}(\BP_*,\BP_*/2)=\F_2[\widetilde{v}_1];\]
see \cite[Th.4.3.2]{greenbook}. The class $\widetilde{v}_1^4$ here is a permanent cycle for degree reasons. Indeed, the first differential $v_1^4$ could support would be a $d_3$, but the only class in bidegree $(7,3)$ supports a $d_3$ and thus can't be a boundary. Everything in higher filtration vanishes on the $E_4$-page.\\

All that is left is to show that $\widetilde{v}_1^4$ is sent to $\overline{u}^4\in \pi_8\ku/2$. This follows from the same argument made in the proof of \Cref{etaisdetectionmodtau}, which shows shows that $\widetilde{v}_1$ is sent to $\overline{u} \in \pi_2 \ku/2$.
\end{proof}

    Using the notation from the above proof, we obtain an infinite family of elements $\partial\widetilde{v}_1^{4k}\in \pi_{8k-1}\Sph$ as the composite
    \begin{equation}\label{geometricrepresentation}
        \Sph^{8k}\to \Sph^{8k}/2\xrightarrow{v_1^{4k}}\Sph/2\xrightarrow{\partial} \Sph^1.
    \end{equation}
    These elements are nonzero as they project to nonzero classes mod $\tau$ in filtration 1; this is because the boundary map
    \[\F_2[\widetilde{v_1}]=\pi_{*,0}\nu\Sph/2,\tau\to\pi_{*-1,1}\nu\Sph/\tau\]
    is an injection when $*>0$, as $\pi_{*,0}\1=\Ext_{\BP_*\BP}^{0,*}(\BP_*,\BP_*)=0$ for $*>0$. The class $\partial\widetilde{v}_1^{4k}$ is then detected by a nonzero class mod $\tau$, and it cannot be hit by a differential because it is in filtration 1. The work of Adams \cite{adamsjx} and Quillen \cite{quillenadams} show this $2$-torsion element is divisible by $2^{\ord_2(k)+3}$. Let us write $\al_k$ for any element in $\pi_{8k-1}\Sph$ such that $2^{\ord_2(k)+3}\al_k=\partial\widetilde{v}_1^{4k}$.\\
    
The odd primary case is simpler. Let $p$ be an odd prime, and write $q=2p-2$ and $v$ for the generator of $\pi_q\ell/p$.
    
\begin{theorem}[{\cite{adamsjx}}]
    For an odd prime $p$, there is a map of spectra $v_1\colon \Sph/p[q]\to \Sph/p$ inducing multiplication by $v$ on $\ell$-homology.
\end{theorem}

\begin{proof}
    The proof is similar to \Cref{adamsperiodicprimetwo}, however, many steps are simplified. Indeed, the cofibre sequence $\Sph \xrightarrow{p} \Sph\to \Sph/p$ can be used to show that $\pi_{2p-2}\Sph/p\simeq \Z/p\Z$. We write $\widetilde{v}_1$ for a generator of this group. To see this class is sent to $v\in \pi_q \ell$, we make the same arguments with the ANSS for $\Sph$ using the left and right units.
\end{proof}


\subsection{The Hurewicz image of $\ko$}\label{dectionofkosubsection}

\begin{theorem}[{\cite{adamsjx}}]\label{detectionofko}
	The unit $f\colon \Sph\to\ko$ is surjective on $\pi_{8k+r}$ for all integers $k$ and $r=1,2$.
\end{theorem}

We will use synthetic spectra in the following proof, both to act as a warm-up for \Cref{detectionresultsforj} and to highlight that although such results can also be obtained using the classical ASS or ANSS for $\ko$, we have still avoided calculations of the comodules $H_\ast(\ko;\F_2)$ and $\BP_\ast \ko$, as well as their associated Ext-groups.

\begin{proof}
There is nothing to show for $k<0$ and $k=0$ we are also done by \Cref{gradedringofktheory}. Let us then set $k>0$ and $r=1$ and consider the following commutative diagram of $\BP$-synthetic homotopy groups
\begin{equation}\label{hurewiczforko}\begin{tikzcd}
    {\pi_{8k+2,0}(\1/\tau)}\ar[r]\ar[d]    &   {\pi_{8k+2,0}(\1/2,\tau)}\ar[r, "\partial"]\ar[d]    &   {\pi_{8k+1,1}(\nu \1/\tau)}\ar[d]    \\
    {\pi_{8k+2,0}(\nu\ko/\tau)}\ar[r]\ar[d]    &   {\pi_{8k+2,0}((\nu\ko)/2,\tau)}\ar[r, "\partial"]\ar[d]    &   {\pi_{8k+1,1}(\nu\ko/\tau})\ar[d]    \\
    {\pi_{8k+2,0}(\nu \ku/\tau)}\ar[r]    &   {\pi_{8k+2,0}((\nu\ku)/2,\tau)}\ar[r, "\partial"]    &   {\pi_{8k+1,1}(\nu \ku/\tau)=0}
\end{tikzcd}\end{equation}
from the maps of $\E_\infty$-rings $\Sph\to \ko\to\ku$ tensored with the cone of $\tau$ and the defining cofibre sequence for $\Sph/2$.\\

Let us write $\widetilde{v}_1\in \pi_2\Sph/2\simeq \Z/4\Z$ for a generator, a lift of $\eta$ through the boundary map $\pi_2 \Sph/2 \to \pi_1 \Sph$, see \Cref{anssmooreattwopicture}, which is detected in $\pi_{2,0}$ of $\1/(2,\tau)$ by a generator which we also denote by $\widetilde{v}_1$. Let us write $x\in \pi_{8k+2}\Sph/2$ for the element $v_1^{4k}\widetilde{v}_1$, and $\mu=\partial(x)\in \pi_{8k+1}\Sph$ for the image of $x$ under the boundary map for the cofibre sequence defining $\Sph/2$. From the proof of \Cref{adamsperiodicprimetwo}, we see that the class $x$ is detected by an element $v_1^{4k+1}$ inside $\pi_{8k+2,0}\1/(2,\tau)$.\\

We claim that the $f(\mu)$ is nonzero in the homotopy groups of $\ko$. To see this, notice that the class $v_1^{4k+1}$ defined above is nonzero, as its image in $\pi_{8k+2,0}(\nu\ku)/(2,\tau)$ is precisely the class $\bar{u}^{4k+1}$, where $\bar{u}$ is the only class on the $E_2$-page of the spectral sequence associated to $\nu\ku/2$ which could detect the generator of $\pi_2 \ku/2$. Indeed, this follows from the fact that the generator $\widetilde{v}_1$ of $\pi_2\Sph/2\simeq \Z/4\Z$ is sent to $u$ from the proof of \Cref{etaisdetectionmodtau}. As the unit $\1\to \nu\ku$ factors through $f$, we see that $f(v_1^{4k+1})$ is also nonzero. As $\pi_{8k+2,0}(\nu\ko/\tau)=0$, see \Cref{manssjattwopictureetwo}, we see that $\partial f(v_1^{4k+1})=f(\partial(v_1^{4k+1}))$ must be nonzero in $\pi_{8k+1,1}(\nu\ko)/\tau$. We now claim that this nonzero element on the $E_2$-page of the $\sigma$-SS from $\nu\ko$ is a permanent cycle. This follows from the fact that we have an explicit geometric representative for this class, $f(\mu)=f(\partial(x))$. Therefore, $f(\mu)$ is nonzero.\\

The case for $r=2$ follows from the above argument by multiplication by $\eta$.
\end{proof}

Notice that our arguments above actually show that $\mu\in \pi_{8k+1}\Sph$ and $\mu\eta\in \pi_{8k+2}\Sph$ are also both nonzero.


\section{Detection results for $\j$}\label{detectionresultsforj}

Extending the ideas used in the proof of \Cref{detectionofko} from $\nu\ko$ to $\j_E$, we will prove \Cref{maintheorem} that the unit map $\Sph\to \j$ induces a surjection on homotopy groups.


\subsection{Hurewicz image of $\j$ at $p=2$}

Implicitly complete all spectra at the prime $2$.

\begin{prop}\label{lowdimensionaldetection}
The unit map $\Sph\to \j$ induces an isomorphism on $\pi_k$ for $k\leq 3$.
\end{prop}

\begin{proof}
    The orders of the homotopy groups in question match, so it suffices to show that $\j$ detects both $\eta\in \pi_1\Sph$ and $\nu\in \pi_3\Sph$. By assumption, we know that $\eta$ is detected in $\pi_1\ko$, and since the unit $\Sph\to \ko$ factors through $\j$, we see $\j$ also detects $\eta$. For $\nu$, we look at $\j_\BP$; see \Cref{manssjattwopictureeinfty}. In particular, the $2$-extension from $2h_2$ to $h_1^3$, comes from the fact that we know $\pi_3 \j\simeq \Z/8\Z$. As $\eta^3$ lies in $\pi_3 \j$ and $\eta^3=4\nu$ in $\pi_3\Sph$, we see that the image of $\nu$ on $\pi_3 \j$ must be nonzero.
\end{proof}

\begin{prop}\label{harddetectionattwo}
The unit map $\Sph\to \j$ induces a surjection on $\pi_{8k-1}$ for all integers $k$.
\end{prop}

\begin{proof}
The case of $k\leq 0$ is uninteresting, so let $k> 0$. Recall that $\al_k$ is any element in $\pi_{8k-1}\Sph$ such that $2^{\ord_2(k)+3}\al_k=\partial\widetilde{v}_1^{4k}$ using the notation from the proof of \Cref{adamsperiodicprimetwo}. As the $2$-power-torsion of these classes $\al_k$ is precisely the size of the cyclic group $\pi_{8k-1}\j$, which we obtained purely from the action of $\psi^3-1$ on the homotopy groups of $\ko$, the surjectivity of $\Sph\to \j$ will follow if the image of $\partial\widetilde{v}_1^{4k}$ are nonzero in $\pi_\ast \j$.\\

Consider the following commutative diagram of synthetic homotopy groups, induced from the maps of synthetic spectra $\1\to \j_\BP\to \ku$ tensored with the cofibre sequence defining $\1/2$ and $C(\tau)$:

    \[\begin{tikzcd}
    {\pi_{8k,0}\1/\tau}\ar[r]\ar[d]    &   {\pi_{8k,0}\1/(2,\tau)}\ar[r, "\partial"]\ar[d]    &   {\pi_{8k-1,1}\1/\tau}\ar[d]    \\
    {\pi_{8k,0}(\j_\BP/\tau)}\ar[r]\ar[d]    &   {\pi_{8k,0}( \j_{\BP}/2,\tau)}\ar[r, "\partial"]\ar[d]    &   {\pi_{8k-1,1}(\j_\BP/\tau)}\ar[d]    \\
    {\pi_{8k,0}(\nu \ku/\tau)}\ar[r]    &   {\pi_{8k,0}((\nu\ku)/2,\tau)}\ar[r, "\partial"]    &   {\pi_{8k-1,1}(\nu \ku/\tau)=0}
    \end{tikzcd}\]

We claim, that the image of $\partial\widetilde{v}_1^{4k}$ inside $\pi_{8k-1,1}(\j_\BP/2,\tau)$ is nonzero. Indeed, by construction, $\widetilde{v}_1^{4k}$ is nonzero in $\pi_{8k,0}((\nu\ku)/2,\tau)$, hence this class is also nonzero in $\pi_{8k,0}(\j_\BP/2,\tau)$. It then follows that $\partial \widetilde{v}_1^{4k}$ is nonzero in $\pi_{8k-1,1}(\j_\BP/2,\tau)$ as we can easily calculate $\pi_{8k,0}(\j_\BP/\tau)=0$ by construction; see \Cref{evenanssjcalculation} and \Cref{manssjattwopictureeinfty}. This shows there is a nonzero class in the $E_2$-page of a modified ANSS for $\j$ detecting the desired element from $\pi_{8k-1}\Sph$, and again, this class is a permanent cycle as it has an explicit geometric representation; simply post-compose \Cref{geometricrepresentation} with the unit $\Sph\to \j$.
\end{proof}

\begin{prop}\label{restofdifficultfamily}
The unit map $\Sph\to \j$ induces a surjection on $\pi_{8k}$ for all $k\in\Z$.
\end{prop}

\begin{proof}
For $k\leq 0$, this is uninteresting or covered by \Cref{lowdimensionaldetection}, so let $k\geq 1$. The class generating $\pi_{8k-1}\j$ supports multiplication by $\eta$. Therefore, the generator of $\pi_{8k}\j\simeq \Z/2\Z$ is hit by the element $\eta\alpha_k$ from the sphere---again, this also proves that $\eta\alpha_k$ is nonzero.
\end{proof}

In degrees $8k+r$ for $r=1,2$, we can deduce the surjectivity of $\Sph\to \j$ using $\ko$.

\begin{prop}\label{detectionattwoonetwo}
The unit map $\Sph\to \j$ induces a surjection on $\pi_{8k+r}$ for all $k\in \Z$ and $r=1,2$.
\end{prop}

\begin{figure}[h]\begin{center}
\makebox[\textwidth]{\includegraphics[trim={2.5cm 15cm 2.5cm 4cm},clip,page = 3, scale = 0.8]{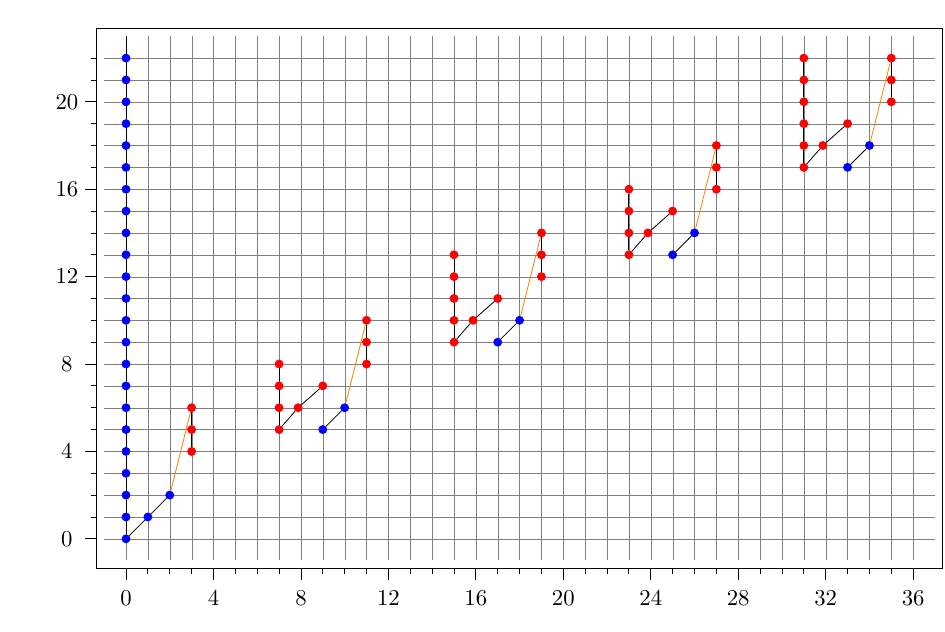}}
\caption{\label{massjattwopictureeinfty} MASS for $\j$ at the prime $2$ (from $\j_{\F_2}$) in the range $0\leq s\leq 36$ on the $E_\infty$-page.}
\end{center}\end{figure}

\begin{figure}[h]\begin{center}
\makebox[\textwidth]{\includegraphics[trim={2.5cm 20.5cm 2.5cm 4cm},clip,page = 2, scale = 0.8]{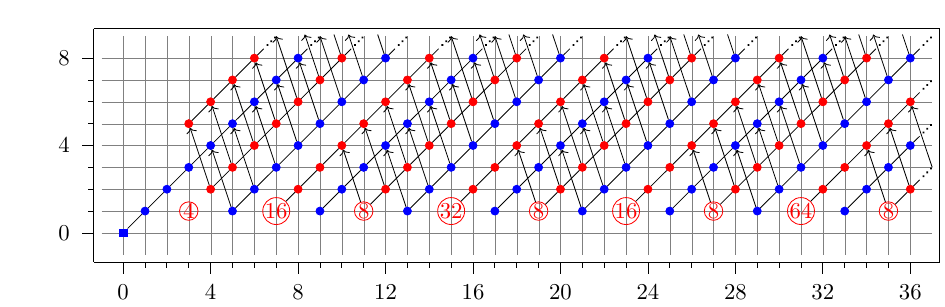}}
\caption{\label{manssjattwopictureetwo} MANSS for $\j$ at the prime $2$ (from $\j_{\BP}$) in the range $0\leq s\leq 36$ with all differentials.}
\end{center}\end{figure}

\begin{figure}[h]\begin{center}
\makebox[\textwidth]{\includegraphics[trim={2.5cm 21.5cm 2.5cm 4cm},clip,page = 3, scale = 0.8]{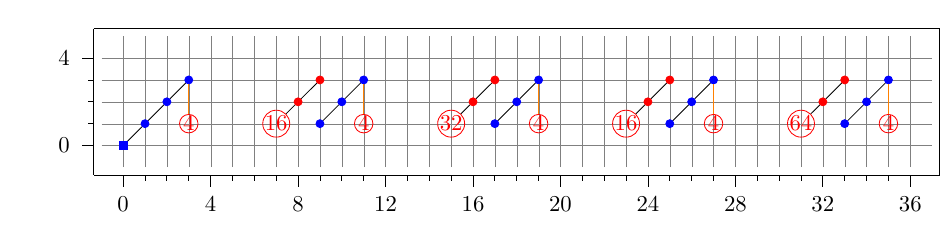}}
\caption{\label{manssjattwopictureeinfty} MANSS for $\j$ at the prime $2$ (from $\j_{\BP}$) in the range $0\leq s\leq 36$ on the $E_\infty$-page.}
\end{center}\end{figure}

\begin{proof}
By construction, the unit map $\Sph\to \ko$ factors through $\j$. By \Cref{detectionofko}, we see that the blue classes in our modified ASS for $\j$, see \Cref{massjattwopictureeinfty}, are hit by elements from the sphere. Combining this with \Cref{restofdifficultfamily}, we know the red class in $\pi_{8k+1}\j$ is hit by $\eta\alpha_k$, which gives us our desired surjectivity.
\end{proof}

\begin{prop}\label{detectionattwothree}
The unit map $\Sph\to \j$ induces a surjection on $\pi_{8k+3}$ for all $k\in \Z$.
\end{prop}

\begin{proof}
The orange extensions by $\eta$ in our modified ASS for $\j$ (\Cref{massjattwopictureeinfty}) combined with \Cref{lowdimensionaldetection} and \Cref{detectionattwoonetwo} show that $\j$ detects $\eta^3\al_k$. This implies that $\j$ also detects $\nu\alpha_k$, as we have $4\nu=\eta^3$ in the homotopy groups of both $\Sph$ and $\j$.
\end{proof}

\begin{remark}\label{periodicityoperatorfailsnot}
Much of the surjectivity of $\Sph\to \j$ is concluded from the classical $\F_2$-ASS for $\j$ using Adams periodicity operator $P=\langle \sigma,16,  -\rangle$. This approach can also be applied here, using the synthetic Toda bracket $\langle  h_3, h_0^4,-\rangle$, where $h_3$ is considered as an element in $\pi_{7,1}\1$, rather than in $\pi_{7,5}\j_{\F_2}$. This is the reason we named the classes in \Cref{evenassjcalculation} with $P^k(-)$ decorations. This operator also produces an alternative proof to many of the statements leading to our proof of \Cref{maintheorem}. We have chosen to give the $\BP$-synthetic proofs here instead as they are simple arguments with long exact sequences.
\end{remark}


\subsection{Hurewicz image of $\j$ at odd primes}

The odd primary case can be proven in one go.

\begin{prop}\label{surjectivityatoddprimes}
    For an odd prime $p$, the unit map $\Sph\to \j$ induces a surjection on homotopy groups.
\end{prop}

\begin{proof}
The proof is similar to \Cref{harddetectionattwo}, so we omit some details. First, we note that $\pi_n \j=0$ unless $n=0$ or $n=kq-1$ for some $k\geq 1$. The key step is now to show that $\partial v_1^k$ is detected in $\pi_{kq-1}\j$. To see this, we use the commutative diagram of synthetic homotopy groups
    \[\begin{tikzcd}
    {\pi_{kq,0}\1/\tau}\ar[r]\ar[d]    &   {\pi_{kq,0}\1/(p,\tau)}\ar[r, "\partial"]\ar[d]    &   {\pi_{kq-1,1}\1/\tau}\ar[d]    \\
    {\pi_{kq,0}(\j_\BP/\tau)}\ar[r]\ar[d]    &   {\pi_{kq,0}( \j_{\BP}/p,\tau)}\ar[r, "\partial"]\ar[d]    &   {\pi_{kq-1,1}(\j_\BP/\tau)}\ar[d]    \\
    {\pi_{kq,0}(\nu \ku/\tau)}\ar[r]    &   {\pi_{kq,0}((\nu\ku)/p,\tau)}\ar[r, "\partial"]    &   {\pi_{kq-1,1}(\nu \ku/\tau)=0}
    \end{tikzcd}\]
and the fact that $\pi_{a,b}\j_\BP/\tau$ is supported in bidegrees $a=kq-1$ and $b=1$ by construction; see \Cref{oddjcalculation}.
\end{proof}

Combining these results yields our main theorem.

\begin{proof}[Proof of \Cref{maintheorem}]
At odd primes this is \Cref{surjectivityatoddprimes}, and at the prime $2$ combine \Cref{harddetectionattwo,restofdifficultfamily,detectionattwoonetwo,detectionattwothree} with the fact that $\pi_{8k+r}\j=0$ for $r\in\{4,5,6\}$.
\end{proof}


\subsection{Split surjectivity on filtered homotopy groups}

These detection statements imply a split surjectivity statement for our modified ANSS for $\j$.

\begin{cor}\label{splittingsurjectiononeinfty}
Let $p$ be a prime. The unit map $\1\to \j_\BP$ induces a split surjection on $E_\infty$-pages of $\sigma$-SSs.
\end{cor}

\begin{proof}
By \Cref{maintheorem}, we see the associated map is surjective. Combined with the calculation of the $\alpha$-family in the ANSS for $\Sph$, see \cite[Th.5.3.7]{greenbook}, we see that the $E_\infty$-page of the $\sigma$-SS for $\j_{\BP}$ detects exactly the $\alpha$-family in $\Sph$, a split summand of the $E_\infty$-page of the ANSS for $\Sph$.
\end{proof}

\begin{remark}\label{failiureofsurjectivity}
Although $\1\to \j_\BP$ induces a surjection on synthetic homotopy groups, this map is not a split surjection. Indeed, notice that $\pi_{3,5}\1=0$, however, $\pi_{3,5}\j\neq 0$ as it contains $\tau$-torsion classes; this is evident from \Cref{evenanssjcalculation} and \Cref{manssjattwopictureetwo}. One might try to clear out some of these $\tau$-torsion classes by using something such as $\tau_{\geq 4}^L \ko$ where $L$ is defined by the equation $x=y$, however, then the line of slope $1$ of red classes supporting the differentials which hit the line we just cut out would survive and not be hit by synthetic unit.\\

On the other hand, it is easy to see that $\1\to \j_{\F_2}$ cannot be split surjective on synthetic homotopy groups for degree reasons.
\end{remark}

\begin{remark}\label{splittingusingJ}
In \cite{adamsjx}, Adams shows that the image of the $J$-homomorphism is a split summand of $\pi_\ast\Sph$. In fact, more is shown, as it is further shown that the $\mu$-family is also a split summand of $\pi_\ast \Sph$. Adams uses the \emph{$e$-invariant} to construct a splitting of the image of the $J$-homomorphism and the $\mu$-family, and according to \cite[\textsection4.1]{kochman}, Mahowald first defined the connective image-of-$J$ spectrum $\j$ in \cite{mahowaldimageofjspectrum} such that the unit map $\Sph\to \j$ lifts this $e$-invariant. To be explicit about Adams splitting now, let us write $i\colon \pi_\ast\j\to \pi_\ast\Sph$ for the Adams map of graded abelian groups which splits the surjective unit map $\pi_\ast \Sph\to \pi_\ast \j$.
\end{remark}

From these considerations, \Cref{filtereredsplitness} follows rather formally from an algebraic fact. Recall that for a filtered abelian group
\[G=G_0\supseteq G_1\supseteq\cdots\]
with $\bigcap G_i=\{0\}$ there is a function $\pi_G\colon G\to \gr G$ from the filtered abelian group to its associated graded abelian groups, which sends $0\mapsto 0$ and for $0\neq x\in G_i\setminus G_{i+1}$, $x$ is sent to $[x]\in G_i/G_{i+1}$. For two such filtered abelian groups $H,G$, given a map of graded abelian groups $j\colon\gr H\to \gr G$, we say a homomorphism of abelian groups $i\colon H\to G$ is a \emph{lift} of $j$ if the following diagram of sets commutes:

\begin{equation}\label{filtereddiagram}
\begin{tikzcd}
H\arrow[d,"\pi_H"]\arrow[r,"i"]&G\arrow[d,"\pi_G"]\\
\gr H\arrow[r,"j"]&\gr G
\end{tikzcd}
\end{equation}

\begin{lemma}\label{splittingoffiltrationlemma}
With notation as above, if $i$ is a lift of $j$, then $i$ is a map of filtered abelian groups and $\gr (i)=j$.
\end{lemma}

\begin{proof}
It suffices to show that if $0\neq x\in H_i\setminus H_{i+1}$, then $i(x)\in G_i$. This follows immediately from the definition of $\pi_G$, since the diagram (\ref{filtereddiagram}) commutes.
\end{proof}

\begin{proof}[Proof of \Cref{filtereredsplitness}]
The unit map $\pi_\ast\Sph\to \pi_\ast\j$ is a morphism of filtered abelian groups where $\pi_\ast\Sph$ and $\pi_\ast\j$ have filtrations defined by $\sigma$-SS of the $\BP$-synthetic spectra $\1$ and $\j_{\BP}$. It suffices to show that the classical splitting map $i$ of \Cref{splittingusingJ} is a morphism of filtered abelian groups. To this end, we apply \Cref{splittingoffiltrationlemma}, where $i$ to be the splitting of \Cref{splittingusingJ} and $j$ as the splitting on $E_\infty$-pages from \Cref{splittingsurjectiononeinfty}. One checks directly that the diagram (\ref{filtereddiagram}) commutes: if $i(x)=y$, then $j$ sends the element on $E_\infty$ detecting $x$ to the element detecting $y$.
\end{proof}

\subsection{The $\alpha$-family and inverting $\eta$} 

In \cite{invertingeta}, Andrews and Miller calculated the bigraded homotopy groups of the $\mathbf{C}$-motivic sphere spectrum after inverting the motivic Hopf map $\eta$. They proceed by calculating the localized Ext groups
\[h_1^{-1}\Ext_{\BP_*\BP}(\BP_*,\BP_*)\]
and show that these are given as $\F_2[h_1^{\pm},\alpha_3,\alpha_4]/(\alpha_4^2)$, where $\alpha_n$ is the $n$-th element of the $\alpha$ family, the same $\alpha_n$'s appearing in our \Cref{evenanssjcalculation}. A crucial point of their calculation is the fact that, in $\Ext_{\BP_*\BP}(\BP_*,\BP_*)$, the elements $h_1^k\alpha_n$ are nonzero, for all $k\ge0$ and $n\neq 2$. Andrews--Miller prove this by way of the algebraic Novikov spectral sequence, but this particular fact goes back to Miller--Ravenel--Wilson \cite[Cor.4.23]{mrw}.\\ 

Via the isomorphism
\[\pi_{*,*}\1/\tau\cong \Ext^{*,*}_{\BP_*\BP}(\BP_*,\BP_*)\]
in $\Syn_\BP$, Andrews--Miller's calculation can thus be interpreted as giving a computation of $\pi_{*,*}\1[h_1^{-1}]$ in $\Syn_\BP$. This relationship between $\mathbb C$-motivic spectra and $\Syn_\BP$ is a shadow of a theorem of Pstragowski, which gives an equivalence of $\infty$-categories between (an even variant) of $\Syn_\BP$ and the category of cellular $\mathbb C$-motivic spectra, after completion at a prime \cite{syntheticspectra}. We can bypass the Ext calculations of Andrews--Miller and Miller--Ravenel--Wilson and deduce the fact that $h_1^k\alpha_n\neq0$ immediately from our synthetic detection theorem for $\j_\BP$.

\begin{theorem}\label{millerandrewsthm}
In $\Syn_\BP$, we have
\[0\neq h_1^k\alpha_n\in \pi_{*,*}\1\]
for all $k\ge0$ and $n\neq 2$.
\end{theorem}

\begin{proof}
If we had $0=h_1^k\alpha_n\in \pi_{*,*}\1$, then the same would be true in $\pi_{*,*}\j_\BP$, contradicting \Cref{evenanssjcalculation}. This is most easily understood by looking at \Cref{manssjattwopictureetwo}; the $\alpha_n$'s appear there in the 1-line. Aside from $\alpha_2$ in bidegree $(3,1)$, each such class supports an infinite $h_1$-tower.
\end{proof}


\addcontentsline{toc}{section}{References}
\bibliography{main} 
\bibliographystyle{alpha}


\appendix
\section{Spectral sequence diagrams}\label{sssection}

The spectral sequences diagrams in this section use the following conventions:

\begin{itemize}
\item Each dot represents a single copy of $\F_p$. Squares represent a copy of the integers $\Z$.
\item If the spectral sequence is the $\sigma$-SS of a synthetic spectrum $X$, defined from either a fibre sequence $X\to Y\to Z$ or a cofibre sequence $U\to W\to X$, then we write in blue those classes detected by the map preserving degree, so $X\to Y$ or $W\to X$, and in red those classes detected by a boundary map, so $Z[-1]\to X$ or $X\to U[1]$. For example, the image-of-$J$ spectra $\j_E$ and $\J_E$ are defined as fibres, where as the mod $2$ Moore spectrum $\Sph/2$ is defined as a cofibre.
\item Orange lines indicate exotic extensions, meaning extensions between two classes of different colours.
\end{itemize}


\begin{figure}[h]\begin{center}
\makebox[\textwidth]{\includegraphics[trim={2.5cm 14.5cm 2.5cm 4cm},clip,page = 4, scale = 0.8]{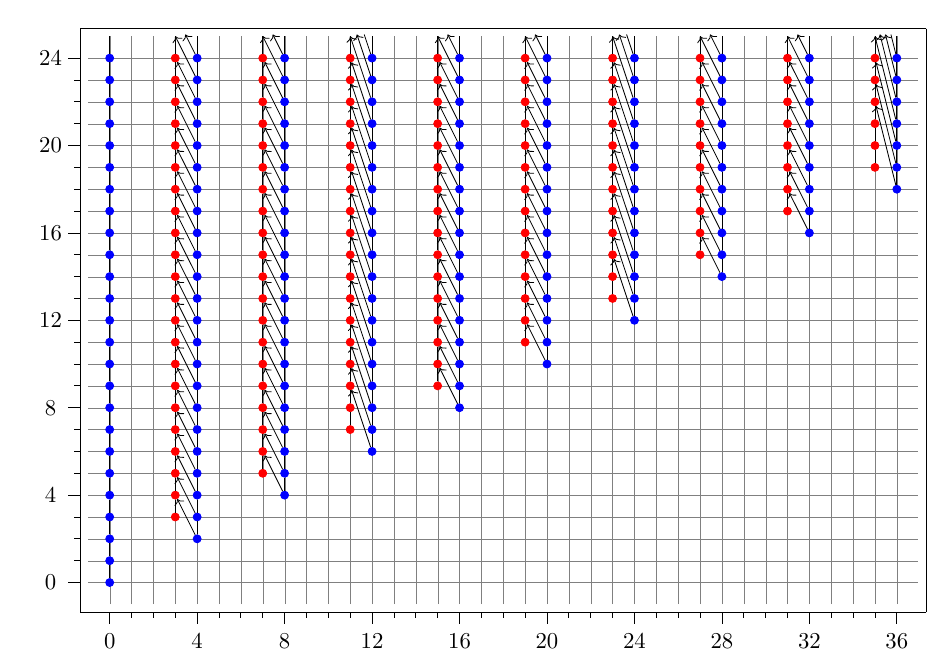}}
\caption{\label{massjatthreepicture} MASS for $\j$ at the prime $3$ (from $\j_{\F_3}$) in the range $0\leq s\leq 36$ with all differentials.}
\end{center}\end{figure}

\begin{figure}[h]\begin{center}
\makebox[\textwidth]{\includegraphics[trim={2.5cm 14.5cm 2.5cm 4cm},clip,page = 5, scale = 0.8]{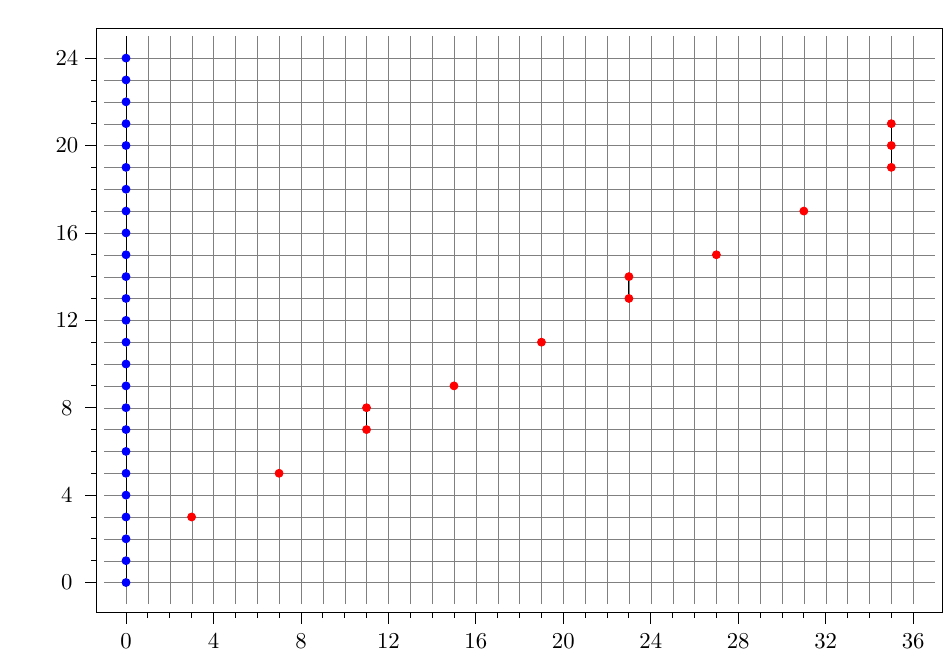}}
\caption{\label{massjatthreepicture} MASS for $\j$ at the prime $3$ (from $\j_{\F_3}$) in the range $0\leq s\leq 36$ on the $E_\infty$-page.}
\end{center}\end{figure}

\begin{figure}[h]\begin{center}
\makebox[\textwidth]{\includegraphics[trim={2.5cm 22.5cm 2.5cm 4cm},clip,page = 4, scale = 0.8]{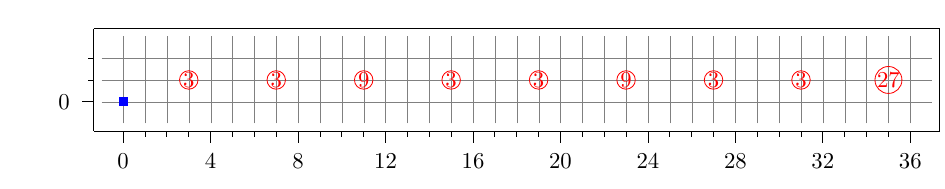}}
\caption{\label{anssjatthreeeinftypicture} MANSS for $\j$ at the prime $3$ (from $\j_\BP$) in the range $0\leq s\leq 36$.}
\end{center}\end{figure}

\begin{figure}[h]\begin{center}
\makebox[\textwidth]{\includegraphics[trim={2.5cm 16cm 2.5cm 4cm},clip,page = 6, scale = 0.8]{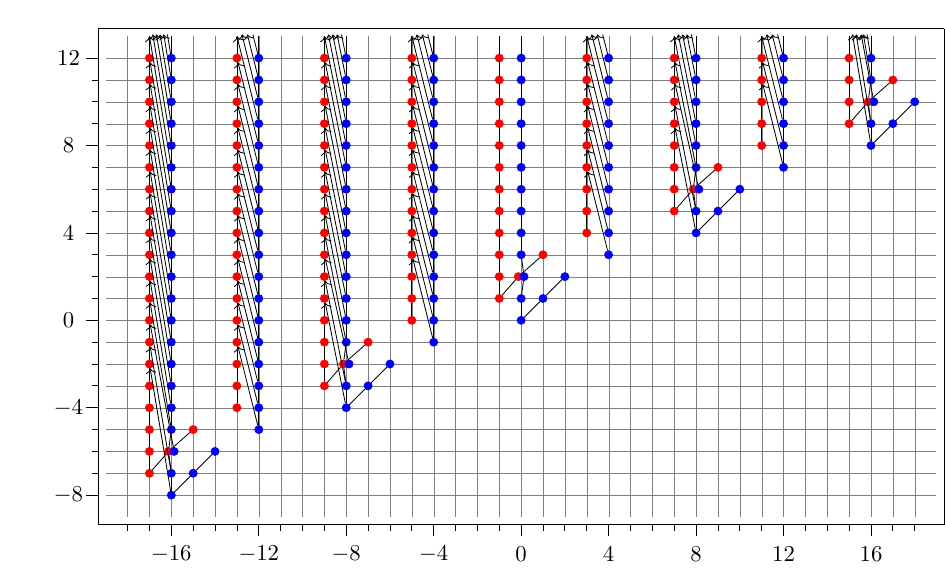}}
\caption{\label{assperiodicjatinftytwoetwopicture} MASS for $\J$ at the prime $2$ (from $\J_{\F_2}$) in the range $-18\leq s\leq 18$ with all differentials.}
\end{center}\end{figure}

\begin{figure}[h]\begin{center}
\makebox[\textwidth]{\includegraphics[trim={2.5cm 16cm 2.5cm 4cm},clip,page = 7, scale = 0.8]{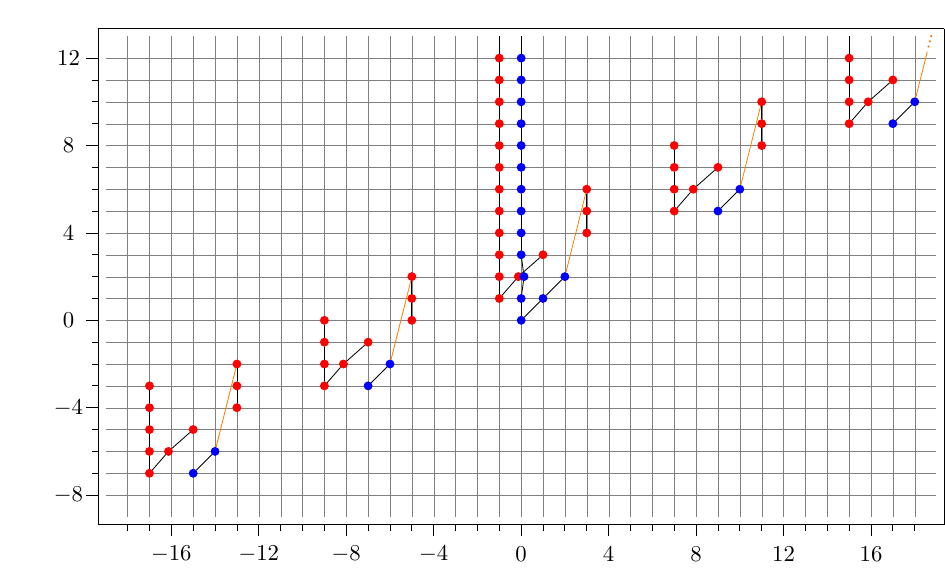}}
\caption{\label{assperiodicjattwoeinftypicture} MASS for $\J$ at the prime $2$ (from $\J_{\F_2}$) in the range $-18\leq s\leq 18$ on the $E_\infty$-page.}
\end{center}\end{figure}

\begin{figure}[h]\begin{center}
\makebox[\textwidth]{\includegraphics[trim={2.5cm 20.5cm 2.5cm 4cm},clip,page = 5, scale = 0.8]{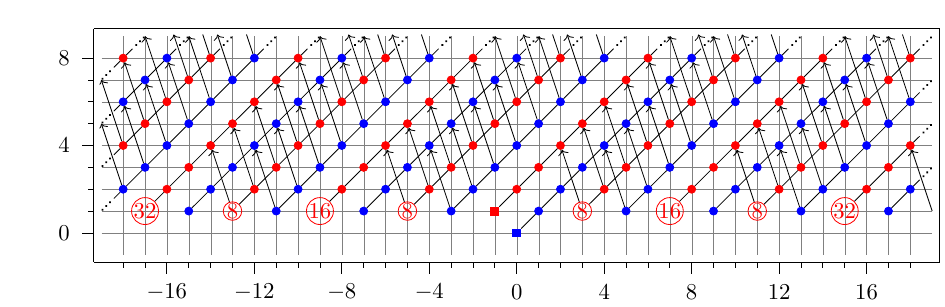}}
\caption{\label{anssperiodicjattwoetwopicture} MANSS for $\J$ at the prime $2$ (from $\J_\BP$) in the range $-18\leq s\leq 18$ with all differentials.}
\end{center}\end{figure}

\begin{figure}[h]\begin{center}
\makebox[\textwidth]{\includegraphics[trim={2.5cm 21.5cm 2.5cm 4cm},clip,page = 6, scale = 0.8]{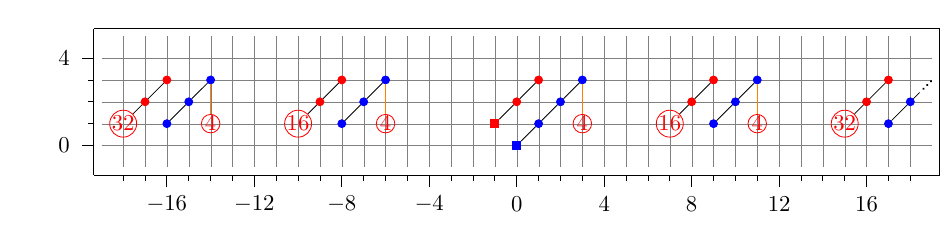}}
\caption{\label{anssperiodicjeinftypicture} MASS for $\J$ at the prime $2$ (from $\J_\BP$) in the range $-18\leq s\leq 18$ on the $E_\infty$-page.}
\end{center}\end{figure}

\begin{figure}[h]\begin{center}
\makebox[\textwidth]{\includegraphics[trim={6cm 21cm 5cm 4cm},clip,page = 8, scale = 0.8]{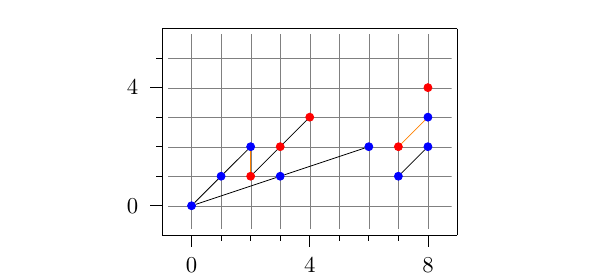}}
\caption{\label{assmooreattwopicture} ASS for $\Sph/2$ at the prime $2$ in the range $0\leq s\leq 8$.}
\end{center}\end{figure}

\begin{figure}[h]\begin{center}
\makebox[\textwidth]{\includegraphics[trim={6cm 20.5cm 5cm 4cm},clip,page = 7, scale = 0.8]{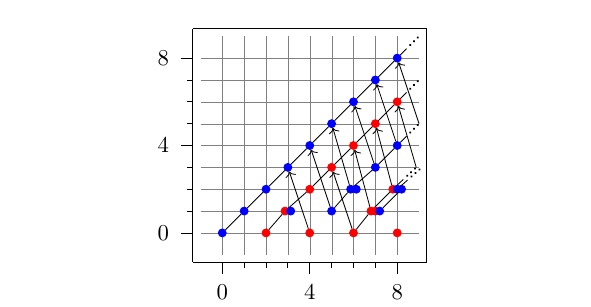}}
\caption{\label{anssmooreattwopicture} ANSS for $\Sph/2$ at the prime $2$ in the range $0\leq s\leq 8$ with all differentials.}
\end{center}\end{figure}

\begin{figure}[h]\begin{center}
\makebox[\textwidth]{\includegraphics[trim={6cm 22cm 5cm 4cm},clip,page = 8, scale = 0.8]{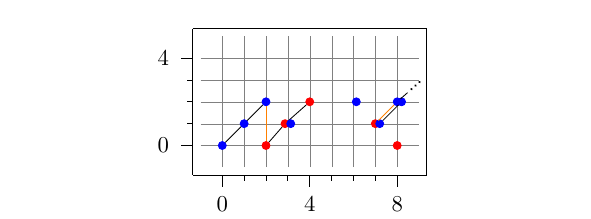}}
\caption{\label{anssmooreattwopictureeinfty} ANSS for $\Sph/2$ at the prime $2$ in the range $0\leq s\leq 8$ on the $E_\infty$-page.}
\end{center}\end{figure}

\end{document}